\documentclass[a4paper,12pt]{article}

\pdfoutput=1

\usepackage[cmtip,arrow]{xy}
\usepackage{pb-diagram,pb-xy}
\usepackage[utf8]{inputenc}
\usepackage[T1]{fontenc}
\usepackage[english]{babel}
\usepackage{amsmath,amssymb,wasysym}
\usepackage{url}
\usepackage{enumerate}
\usepackage{pdfpages}

\title{Etienne Bézout on elimination theory
\footnote{ We wish to thank Christian Houzel and Roshdi Rashed for their advice.}}
\author{Erwan Penchèvre}
\date{}

\begin{document}
\maketitle

Bézout's name is attached to his famous theorem. Bézout's Theorem states that the degree of the eliminand of a system a $n$ algebraic equations in $n$ unknowns, when each of the equations is generic of its degree, is the product of the degrees of the equations. The eliminand is, in the terms of XIXth century algebra\footnote{\textit{cf.} \cite{netto1900} vol. 2.}, an equation of smallest degree resulting from the elimination of $(n-1)$ unknowns.

Bézout demonstrates his theorem in 1779 in a treatise entitled \textit{Théorie générale des équations algébriques}. In this text, he does not only demonstrate the theorem for $n>2$ for generic equations, but he also builds a classification of equations that allows a better bound on the degree of the eliminand when the equations are not generic. This part of his work is difficult: it appears incomplete and has been seldom studied. In this article, we shall give a brief history of his theorem, and give a complete justification of the difficult part of Bézout's treatise.

\section{The idea of Bézout's theorem for $n>2$}

The theorem for $n=2$ was known long before Bézout. Although the modern mind is inclined to think of the theorem for $n>2$ as a natural generalization of the case $n=2$, a mathematician rarely formulates a conjecture before having any clue or hope about its truth. Thus, it is not before the second half of XVIIIth century that one finds a clear statement that the degree of the eliminand should be the product of the degrees even when $n>2$.

Lagrange, in his famous 1770-1771 memoir \textit{Réflexions sur la résolution algébrique des équations} \cite{lagrange1770}, proves Bézout's theorem for several particular systems of more than two equations, by studying the functions of the roots remaining invariant through some permutations. In the same year 1770, Waring enunciates the theorem for more than two equations in his \textit{Meditationes algebraicae} \cite{waring1770}, without demonstration. Up to our knowledge, these are the first occurences of Bézout's theorem for $n>2$. 

Bézout probably knew those works of Lagrange and Waring. He is directly concerned by Lagrange's memoir, where Lagrange nominally criticizes the algebraical methods of resolution of equations in one unknown, that Bézout had conceived in the 1760's. Waring says, in the preface to the second edition of his \textit{Meditationes algebraicae}, having sent a copy of its first edition, as soon as 1770, to some scholars, including Lagrange and Bézout\footnote{\textit{Cf.} p.~\textit{xxi} de l'édition de 1782 des \textit{Meditationes algebraicae}}. Thus Bézout's theorem was already in the mind of those three scholars as soon as 1770.\footnote{It is to be noted, that Lagrange and Waring where also among the first readers of Bézout's treatise in 1779. They will, shortly after, give an answer, Lagrange in his correspondance with Bézout, and Waring in the second edition of his \textit{Meditationes algebraicae}.}

On the contrary, Bézout was not yet aware of the formula of the product of the degrees in 1765. His works \cite{bezout1762,bezout1765} of the years 1762-1765 about the resolution of algebraic equations show several examples
of systems, of more than two equations, where the method designed by him in 1764 leads him to a final equation of a degree much higher than the product of the degrees of the initial equations, because of a superfluous factor. The discovery of Bézout's theorem for $n>2$, still as a conjecture, is thus clearly circumscribed in the years 1765-1770.

\section{Bézout's method of elimination and the superfluous factors}

In elimination theory, early works for $n=2$ already show two different and complementary methods (\cite{euler1748a,cramer1750,bezout1764}). One of them relies upon symetrical functions of the roots. This method was used by Poisson to give an alternative demonstration of Bézout's Theorem in 1802. As we have chosen to concentrate on Bézout's path, we won't describe this method in this article.\footnote{Although for $n=2$, \textit{cf.} footnote \ref{Poisson} p. \pageref{Poisson}. Poisson knew of Bézout's work; but he ascribes to it a lack of rigour, thus justifying his own recourse to a different method. \textit{Cf.} \cite{poisson1802}, p.~199; and \cite{poisson1802}, p.~203. One should understand this judgement after section \ref{Demonstration2} below.}

The other method, the one used by Bézout, is a straightforward generalization\footnote{\label{MethodeNewton}This other method was already used for systems of equations of higher degrees by Newton (\emph{cf.} \cite{newton1972} p.~584-594) and before Newton (\emph{cf.} \cite{penchevre2004a}). Bézout sometimes ascribes this method to Newton.} of the principle of substitution used to eliminate unknowns in systems of linear equations ; this principle is still taught today in high school. 

This method does not dictate the order in which to eliminate unknowns and powers of the unknowns. When Bézout uses this method in 1764, for $n>2$, he eliminates the unknowns \textit{one after the other}. This necessarily leads to a superfluous factor increasing the degree of the final equation far above the product of the degrees of the initial equations. This difficulty is easily illustrated by the following system of three equations : 
\label{FacteursSuperflus} 
\begin{align}
&-x^2+y^2+z^2-2yz-2x-1=0\\
&z+x+y-1=0\\
&z-x+y+1=0
\end{align}
Eliminating $z$ between $(1)$ and $(2)$, one obtains 
\begin{align}
4y^2+4xy-4x-4y=0
\end{align}
Eliminating $z$ between $(1)$ and $(3)$, one obtains 
\begin{align}
4y^2-4xy-4x+4y=0
\end{align}
Eliminating $4xy$ between $(4)$ and $(5)$, one has $x=y^2$. Substituting it for $x$ in $(4)$, the final equation is
$$4y(y^2-1)=0$$
The root $y=0$ does not correspond to any solution of the system above.\footnote{not even an infinite solution, in $\mathbb{P}^{3}$.} In fact, the true eliminand should be $y^2-1=0$.
Bézout was well aware of the difficulty. In 1764, he says :
\begin{quote}
(...) when, having more than two equations, one eliminates by comparing them two by two; even when each equation resulting from the elimination of one unknown would amount to the precise degree that it should have, it is vain to look for a divisor, in any of these intermediate equations, that would lower the degree of the final equation; none of them has a divisor; only by comparing them will one find an equation having a divisor; but where is the thread that would lead out of the maze?\footnote{\textit{Cf.} \cite{bezout1764}, p. 290; and also \cite{bezout1779}, p. vii.
}
\end{quote}
At the time in 1764, Bézout had not yet found the exit out of the maze.

Fifteen years later, in his 1779 treatise, Bézout gets rid of this iterative order by reformulating his method in terms of a new concept called the ``sum-equation'' :
\label{EquationSomme}
\begin{quote}
We conceive of each given equation as being multiplied by a special polynomial. Adding up all those products together, the result is what we call the \textit{sum-equation}. This sum-equation will become the final equation through the vanishing of all terms affected by the unknowns to eliminate.
\footnote{\textit{Cf.} \cite{bezout1779}, \S~224.}
\end{quote}
In other words, for a system of $n$ equations with $n$ unknowns\footnote{Throughout our commentary, we shall use upper indices to distinguish equations or polynomials ($f^{(1)}$, $f^{(2)}$, ...), and lower indices to distinguish unknowns or indeterminates ($x_1$, $x_2$, ...). We are thus losely following Bézout's own notations.}, \textit{viz.}
$$\left\lbrace\begin{aligned}
&f^{(1)}=0\\
&\vdots\\
&f^{(n)}=0
\end{aligned}\right.,$$
Bézout postulates that the final equation resulting from the elimination of $(n-1)$ unknowns is an equation of smallest degree, of the form
$$\phi^{(1)}f^{(1)}+...+\phi^{(n)}f^{(n)}=0.$$
The application of the method is thus reduced to the determination of the polynomials $\phi^{(i)}$. First of all, one must find the degree of those polynomials, as well as the degree of this final equation. This is the ``node of the difficulty'' according to Bézout. Once acertained the degree of the final equation, elimination is reduced to the application of the method of undetermined coefficients, and thus, to the resolution of a unique system of \textit{linear} equations. Hence the need for Bézout's theorem predicting the degree of the final equation.

Although this idea of ``sum-equation'' seems a conceptual break-through reminding us of XIXth century theory of ideals, the immediate effect of this evolution is the rather complicated structure of Bézout's treatise ! For didactical reasons maybe, he introduces his concept of ``sum-equation'' only in the second part of his treatise (``livre second''). In the first part of it, his formulation is a compromise with the classical formulation of the principle of substitution. We shall analyse this order of presentation in section \ref{ThreeViews}.

\section{A treatise of ``finite algebraic analysis''}

In the dedication of his \textit{Théorie générale des équations algébriques}, Bézout says that his purpose is to ``perfect a part of Mathematical Sciences, of which all other parts are awaiting what would now further their progress''\footnote{\textit{Cf.} \cite{bezout1779}.}. In the introduction to the treatise, Bézout opposes two branches of the mathematics of his days: ``finite algebraic analysis'' and ``infinitesimal analysis''. The former is the theory of equations. Historically, it comes first ; according to Bézout, infinitesimal analysis has recently drawn all the attention of mathematicians, being more enjoyable, because of its many applications, and also because of the obstacles met with in algebraic analysis. Bézout says :
\begin{quote}
The former itself [infinitesimal analysis] needs the latter to be perfected.\\
(...)\\
The necessity to perfect this part [algebraic analysis] did not escape the notice of even those to whom infinitesimal analysis is most redeemable.\footnote{\textit{Cf.} \cite{bezout1779}, p. ii.}
\end{quote}
In his view, the logical priority of algebraic analysis adds thus to its historical priority. The composition of his treatise is almost entirely algebraic:
\begin{itemize}
\item Bézout only briefly alludes (\S~48) to the geometric interpretation of elimination methods as research of the intersection locus of curves%
\footnote{although he knew well Euler's memoir of 1748, \textit{Démonstration sur le nombre des points où deux lignes des ordres quelconques peuvent se rencontrer}.};
\item he does never make any hypothesis about the existence or the arithmetical nature of the roots of algebraic equations%
\footnote{It seems that the very notion of root does never make appearance anywhere in his demonstrations. The word appears in \S\S~48, 117, 280--284, but never crucially. Bézout knew, of course, that the known methods of algebraic resolution of equations do not apply beyond fourth degree, as Lagrange had explained it exhaustively in his \textit{Réflexions sur la résolution algébrique des équations}. Moreover, at the time, the status of complex numbers and the fundamental theorem of algebra where still problematic.}. 
\end{itemize}
This position of his treatise as specialized research on algebraic analysis is quite singular for his time%
\footnote{H. Sinaceur has commented upon the use of the term ``analysis'' in XVIIIth century (\cite{sinaceur1991}, p.~51): 
\begin{quote}
The term ``analysis'' is a generic concept for the mathematical method rather than a particular branch of it. It was then normal that no clear distinction should exist between algebra and analysis, nor any exclusive specialization. Moreover, the analysis of equations, also called ``algebraic analysis'', could be considered as a part of a whole named ``mathematical analysis''.
\end{quote}
}.

\section{A classification of equations}

In fact, Bézout does not only demonstrate his theorem for \textit{generic} equations. When the equations are not generic, the degree of the eliminand may be less than the product of the degrees of the equations. Bézout progressively studies larger and larger classes of equations, by asking that some coefficients vanish or verify certain conditions. He thus builds a classification that allows a better bound on the degree of the eliminand according to the species of the equations. The case of generic equations is thus encompassed, as a very special case, in a research of gigantic proportion. In this regard, Bézout says: 
\begin{quote}
Whatever idea our readers might have conceived of the scale of the matter that we are about to study, the idea that he will soon get therefrom, will probably surpass it.\footnote{\textit{Cf.} \cite{bezout1779}, \S~52.}
\end{quote}

\paragraph{An example} In \S~62 of his treatise, Bézout proves everything that was known before, in the case $n=2$. For two equations with two unknowns $x_{1}$, $x_{2}$ of the form
$$
\begin{array}{ll}
\sum\limits _{k_{1}\leq a_{1},\ k_{1}+k_{2}\leq t}{A_{k_{1}k_{2}}x_{1}^{k_{1}}x_{2}^{k_{2}}}=0\\
\sum\limits _{k_{1}\leq a_{1}^{\prime},\ k_{1}+k_{2}\leq t^{\prime}}{A_{k_{1}k_{2}}^{\prime}x_{1}^{k_{1}}x_{2}^{k_{2}}}=0\end{array}
$$
where the $A_{k_1k_2}$ and the $A'_{k_1k_2}$ are undetermined coefficients, 
the degree of the final equation resulting from the elimination of $x_2$ is
$D=tt'-(t-a_{1})(t'-a_{1}^{\prime})$. Cramer in 1750, Euler,
then Bézout himself in 1764, had known this result. 
Specifying $a_1=t$, $a_1'=t'$, one obtains the case of two ``complete equations'', \textit{i. e.} generic of their degrees :
$$D=tt'$$ 
In this case, the degree of the eliminand is the product of the degrees of the initial equations.

\paragraph{Orders and species}
What Bézout calls a ``complete polynomial'' is a polynomial, generic of its degree. Non-generic polynomials are called ``incomplete''. Bézout discriminates between several ``orders'' of incomplete polynomials. He thus defines the ``incomplete polynomials of the first order'' as those verifying the following conditions :\footnote{Bézout is either using the term ``polynomial'' or the term ``equation''; here, an equation is always of the form $f=0$ where $f$ is a polynomial.}
\begin{quote}
  \begin{description}
  \item[1\textsuperscript{o}] that the total number of unknowns being $n$, their combinations $n$ by $n$ should be of some degrees, different for each equation
  \item[2\textsuperscript{o}] that their combinations $n-1$ by $n-1$ should be of some degrees, different not only for each equation, but also for each combination
  \item[3\textsuperscript{o}] that their combinations $n-2$ by $n-2$ should be of some degrees, different not only for each equation, but also for each combination;
  \item etc.\footnote{\textit{cf.} \cite{bezout1779}, p. xiii.}
  \end{description}
\end{quote}

Among polynomials of this order, Bézout distinguishes several ``species'' (by the way, the two equations already mentioned in the example above are from the ``first species of incomplete equations''). He says:
\begin{quote}
As it is not possible to attack this problem from the forefront (the problem of incomplete polynomials of first order), I took it in the inverse order, first supposing the absence of the highest degrees of the combinations one by one, then the absence of those and of the highest degrees of the combinations two by two, etc., and also supposing some restrictive conditions in order to facilitate the intelligence of the method (...)
\end{quote}
We shall soon describe the ``restrictive conditions'' alluded to.

Bézout's symbolic notations for incomplete polynomials are such:
\begin{align*}
&(u^{a}...n)^{t} & \text{(\textit{cf.} {\S}57--67)}\\
&\lbrack(u^{a},x^{\underset{\prime}{a}})^{b},y^{\underset{\prime\prime}{a}}...n\rbrack^{t} & \text{(\textit{cf.} {\S}74--81)}\\
&(\lbrack(u^{a},x^{\underset{\prime}{a}})^{b},(u^{a},y^{\underset{\prime\prime}{a}})^{\underset{\prime}{b}},(x^{\underset{\prime}{a}},y^{\underset{\prime\prime}{a}})^{\underset{\prime\prime}{b}}\rbrack^{c},z^{\underset{\prime\prime\prime}{a}}...n)^{t} & \text{(\textit{cf.} {\S}82--132)}\\
&...
\end{align*}
For example, the second line describes a polynomial that we could write today, in a slightly modernized notation but still keeping with Bézout's unusual underscripts:
$$
\sum\limits_{\begin{array}{l}k\leq a,\underset{'}{k}\leq\underset{'}{a},\underset{''}{k}\leq\underset{''}{a},...,\\
k+\underset{'}{k}\leq b,\\
k+\underset{'}{k}+\underset{''}{k}+...\leq t\end{array}}
A_{k\underset{'}{k}\underset{''}{k}}u^{k}x^{\underset{'}{k}}y^{\underset{''}{k}}...
$$
where $u,x,y,z,...$ are the unknowns.

Finally, in section III of book I, Bézout introduces the second, third, fourth orders, etc. of polynomials, represented by this notation:
$$
(u^{a,\overline{a},\overline{\overline{a}},...}...n)^{t,\overline{t},\overline{\overline{t}},...}
$$
where $a\leq\Bar{a}\leq\Bar{\Bar{a}}\leq...$ and $t\geq\Bar{t}\geq\Bar{\Bar{t}}\geq...$. We could write such a polynomial under the form:
\begin{align*}
& \sum\limits _{k\leq\Bar{\Bar{a}},...,k+\underset{'}{k}+\underset{''}{k}+...\leq\Bar{\Bar{t}}}{A_{k\underset{'}{k}\underset{''}{k}}u^{k}x^{\underset{'}{k}}y^{\underset{''}{k}}...}\\
+ & \sum\limits _{k\leq\Bar{a},...,\Bar{\Bar{t}}<k+\underset{'}{k}+\underset{''}{k}+...\leq\Bar{t}}{A_{k\underset{'}{k}\underset{''}{k}}u^{k}x^{\underset{'}{k}}y^{\underset{''}{k}}...}\\
+ & \sum\limits _{k\leq a,...,\Bar{t}<k+\underset{'}{k}+\underset{''}{k}+...\leq t}{A_{k\underset{'}{k}\underset{''}{k}}u^{k}x^{\underset{'}{k}}y^{\underset{''}{k}}...}\\
+ & ...
\end{align*}

The whole book I of Bézout's treatise is thus progressing in an increasing order of generality. Bézout says, several times, that the equations studied earlier are particular cases of the new forms under study\footnote{\textit{cf.} \cite{bezout1779} \S~64, 81.}.

\paragraph{Four important cases}
We don't want to describe in full generality all the cases studied by Bézout.
We shall concentrate on the four large following classes of equations:
\begin{itemize}
\item complete equations, for all $n$
\item first species of incomplete equations, for all $n$
\item second species of incomplete equations, for all $n$
\item third species of incomplete equations, for $n=3$
\end{itemize}
Only for those four classes, Bézout gives an explicit formula for the degree of the eliminand when each proposed equation is generic within its species. We shall give a detailed summary of Bézout's demonstration for the second species, slightly modernized as regards symbolic notations and algebraic terminology\footnote{\textit{Cf.} sections \ref{Demonstration2} below. Some steps of the demonstration will have to await a complete justification in section \ref{Toric}. As for complete equations and for the first species of incomplete equations, results will be derived from the case of second species ; we also provide a more elementary proof in the appendix when $n=3$. As for the third species of incomplete equations, we shall say more in section \ref{Demonstration3}, with a demonstration in section \ref{Toric}.}.

To describe our notations, let us consider a system of $n$ equations in $n$ unknowns :
$$\left\lbrace\begin{aligned}
f^{(1)}=0\\
f^{(2)}=0\\
\vdots\\
f^{(n)}=0
\end{aligned}\right.$$
where the $f^{(i)}$ are elements of a polynomial ring $C=K[x_1,x_2,...,x_n]$. Bézout himself is using several kinds of indices : upper index means equation number, and lower indices mark unknown quantities\footnote{\textit{cf.} \cite{bezout1779} \S 62.}. 
We shall also use multi-index notations and write $k=(k_1,k_2,...,k_n)$ and $x^k=x_1^{k_1}x_2^{k_2}...x_n^{k_n}$. The support $\text{supp}(f^{(i)})$ of a polynomial is the set of points $k\in\mathbb{Z}^n$ such that the monomial $x^k$ has non-zero coefficient in $f^{(i)}$. The main breakthrough of Bézout is thus to distinguish cases with respect to the convex envelop of $\text{supp}(f^{(i)})$.

Let $t$, $a_1$, $a_2$,..., $a_n$ be integers verifying the following conditions (the ``restrictive conditions'' alluded to, in the quotation above) :
$$\left\lbrace\begin{array}{l}
(\forall i)\ a_i\leq t\\
(\forall i\neq j)\ a_i+a_j>t
\end{array}\right.$$
Let $E_{t,a}$ be the convex set in $\mathbb{Z}^n$ defined by
$$\left\lbrace\begin{array}{l}
0\leq k_1\leq a_1\\
0\leq k_2\leq a_2\\
\vdots\\
0\leq k_n\leq a_n\\
k_1+k_2+...+k_n\leq t
\end{array}\right.$$
For $n=3$, such a convex set is the top polyhedron on figure 4 at the end of this article. For any such $t$ and $a$, we define the sub-vector space of $C$ over $K$ of polynomials with support in $E_{t,a}$:
$$C_{\leq t,a}=\left\lbrace f\in K[x_1,x_2,...,x_n] \mid \text{supp}(f)\subset E_{t,a}\right\rbrace$$
An incomplete equation of the first species is a generic member of $C_{\leq t,a}$.

As for systems of equations, let it be given for any index $i\in\lbrace 1,2,...,n\rbrace$ such a set of integers $t^{(i)}$, $a^{(i)}_1$, $a^{(i)}_2$,..., $a^{(i)}_n$, and $f^{(i)}$ be a generic member of 
$$C_{\leq t^{(i)},a^{(i)}}$$
In other words, $a^{(i)}_j=\deg_j f^{(i)}$ is the degree of $f^{(i)}$ with respect to $x_j$, and $t^{(i)}=\deg f^{(i)}$ is the total degree of $f^{(i)}$. For all $k\neq j$, we require that 
$$\deg_jf^{(i)}+\deg_kf^{(i)}\geq\deg f^{(i)}$$
Finally, the genericity of $f^{(i)}$ actually means that
$$\text{supp}(f^{(i)})=E_{t^{(i)},a^{(i)}}$$ 
and that the non-zero coefficients of $f^{(i)}$ are indeterminates adjoined to a base field. For example, let us say $K$ is purely transcendant over $\mathbb{Q}$ :
$$K=\mathbb{Q}((u^{(i)}_k)_{1\leq i\leq n,\ k\in\text{supp}(f^{(i)})})$$ 
where every indeterminate $u^{(i)}_k$ is the coefficient of $x^k$ in $f^{(i)}$. The equations $f^{(i)}=0$ are what Bézout calls a system of ``incomplete equations of the first species''.

As for the ``second species of incomplete equations''\footnote{\textit{cf.} \cite{bezout1779} \S 74-81.}, let $t$, $a_1$, $a_2$,..., $a_n$, $b$ be non-negative integers satisfying the following restrictive conditions:\label{secondspecies}
$$\left\lbrace\begin{array}{l}
\max(a_1,a_2)\leq b\\
(\forall i\geq 3)\ a_i\leq t\\
a_1+a_2\geq b,\\
(\forall i\neq j,\ \lbrace i,j\rbrace\neq\lbrace 1,2\rbrace)\ a_i+a_j\geq t\\
b\leq t\\
(\forall i\geq 3)\ a_i+b\geq t
\end{array}\right.$$
These integers define a convex set $E_{t,a,b}$ in $\mathbb{Z}^n$:
$$\left\lbrace\begin{array}{l}
0\leq k_1\leq a_1,\quad 0\leq k_2\leq a_2,...,\ 0\leq k_n\leq a_n,\\
k_1+k_2\leq b\\
k_1+k_2+...+k_n\leq t
\end{array}\right.$$
For any such $t,a,b$, we define the sub-vector space $C_{\leq t,a,b}$ of $C$ over $K$ of polynomials with support in $E_{t,a,b}$. An incomplete equation of the second species is a generic member of $C_{\leq t,a,b}$.

As for the ``third species of incomplete equations''\footnote{\textit{cf.} \cite{bezout1779} \S 82-132.}, when $n=3$, let $t$, $a_1$, $a_2$, $a_3$, $b_1$, $b_2$, $b_3$ be non-negative integers satisfying the following restrictive conditions:
$$\left\lbrace\begin{array}{l}
\max(a_1,a_2)\leq b_3,\quad\max(a_1,a_3)\leq b_2,\quad\max(a_2,a_3)\leq b_1\\
a_1+a_2\geq b_3,\quad a_1+a_3\geq b_2,\quad a_2+a_3\geq b_1\\
\max(b_1,b_2,b_3)\leq t\\
\min(a_1+b_1,a_2+b_2,a_3+b_3)\geq t\\
b_1+b_2+b_3\geq 2t
\end{array}\right.$$
These integers define a convex set $E_{t,a,b}$ in $\mathbb{Z}^3$:
$$\left\lbrace\begin{array}{l}
0\leq k_1\leq a_1,\quad 0\leq k_2\leq a_2,\quad 0\leq k_3\leq a_3,\\
k_1+k_2\leq b_3,\quad k_1+k_3\leq b_2,\quad k_2+k_3\leq b_1,\\
k_1+k_2+k_3\leq t
\end{array}\right.$$
An incomplete equation of the third species is a generic polynomial equation with support in any such convex set. Bézout further subdivides the third species into eight ``forms'', according to the algebraic signs of the quantities
$$t-b_1-b_2+a_3,\quad t-b_2-b_3+a_1,\quad t-b_3-b_1+a_2$$
His second, third and seventh forms are exchanged under permutation of the unknowns, as well as his fourth, fifth and eighth forms. The four lower polyhedra on figure 4 at the end of this article represent supports belonging to the first, the third, the fifth and the sixth forms.

\section{Three different views on elimination}\label{ThreeViews}

The early article by Bézout (1764) was based on the following idea. The elimination process is split into a sequence of operations. Each operation consists in multiplying some of the previously obtained equations by suitable polynomials, and then building the sum of the products. This idea was already known. Yet, it is perfected in book II of Bézout's 1779 treatise. There, Bézout suggests to represent directly the final equation resulting from the elimination, as such a sum of products of the initial equations by suitable polynomials. This representation could well seem natural to the modern reader used to the concept of ``ideal'' inherited from end-of-nineteenth-century algebra. Analysing the proofs of different cases of Bézout's theorem such as he wrote them in his treatise, we are going to study three different views along the path leading to Bézout's concept of \textit{sum-equation}.

We shall limit ourselves in this section to the case of three equations with three unknowns:
$$\left\lbrace\begin{aligned}
f^{(1)}=0\\
f^{(2)}=0\\
f^{(3)}=0
\end{aligned}\right.$$
of respective degrees $t^{(1)},t^{(2)},t^{(3)}$. 
Bézout postulates at first the existence of a unique final equation resulting from the elimination of two unknowns (\textit{e.~g.} $x_2$ and $x_3$), of minimal degree $D$. Among the three different views that we are going to expound, the two former ones consist in multiplying $f^{(1)}$ by a ``multiplier-polynomial'' $\phi^{(1)}$ with undetermined coefficients, of degree $(T-t^{(1)})$, where $T$ is a large enough integer, and then making all monomials ``vanish'' in the product $\phi^{(1)}f^{(1)}$, except monomials $1,\ x_1,\ x_1^2,...,\ x_1^D$. Those monomials are building the final equation that we are looking for\footnote{Bézout proceeds in this way for ``incomplete equations of the first species'' for example, \textit{cf.} \cite{bezout1779}, \S~59--67.}. This vanishing could be operated in two steps : first use equations $(2)$ and $(3)$ to make vanish as many monomials as possible, and then use the classical method of undetermined coefficients to do the rest of it. After having used equations $(2)$ and $(3)$, in order to apply the method of undetermined coefficients, there should remain as many vanishable terms (each one of them gives an equation), as undetermined coefficients in $\phi^{(1)}$. But the situation is complicated by the fact that only some of the
undetermined coefficients provided by $\phi^{(1)}$ could, according to Bézout,
serve the purpose of elimination, several of them being ``useless''.
Finally, Bézout ends up with a formula:
$$(\text{number of terms remaining in }\phi^{(1)}f^{(1)})-D=\text{nbr. of useful coefficients in }\phi^{(1)}$$
He would then use this formula to calculate $D$.

As we can see, this
method is quite ambiguous, as long as we don't give a more precise meaning to the word ``useless''. Bézout says that the number of ``useless'' coefficients in $\phi^{(1)}$ is precisely the number of monomials that we could make vanish in $\phi^{(1)}$ using $(2)$ and $(3)$. This is, at best, mysterious, as long as we don't give a precise meaning to all this in terms of dimensions of vector spaces. Bézout does not say anything to clarify this situation, as he reserves the effective calculation for book II.

\paragraph{Substituting monomials}
We have said that there are three different views on elimination in Bézout's treatise, and we have just expounded the common setting of the first two of
them. What differentiates them is the way to count the number of 
monomials that we could make vanish in a given polynomial,
thanks to given equations.
In his \emph{first} proof\footnote{As he does for ``complete equations'' (complete, \textit{i. e.} generic of their degrees), \textit{cf.} \cite{bezout1779} \S~45--48.}, Bézout says that: 
$$\text{nbr. of vanishable terms in }\phi^{(1)}f^{(1)}=\text{nbr. of terms divisible by }x_{2}^{t^{(2)}}\text{ or }x_{3}^{t^{(3)}}$$
This reminds us of Newton's method of elimination of highest powers of the unknown\footnote{\emph{Cf. supra} note \ref{MethodeNewton} p.~\pageref{MethodeNewton}.}, generalized to any number of equations and unknowns. It is remarkable that
in 1764, after having said that Euler's and Cramer's methods of elimination\footnote{\label{Poisson}When Bézout mentions Euler in \cite{bezout1764}, he is surely refering to \cite{euler1748a}. For a quick overview of Euler's method in \cite{euler1748a}, let there be two equations of degrees $m$ and $n$ in $y$:
$$\left\lbrace\begin{array}{l} y^m-Py^{m-1}+Qy^{m-2}-...=0\\ y^n-py^{n-1}+qy^{n-2}-...=0 \end{array}\right.$$
where $P$,$Q$,...,$p$,$q$,... are polynomials in $x$, such that:
$$\begin{array}{l} \deg P+m-1=\deg Q+m-2=...=m\\ \deg p+n-1=\deg q+n-2=...=n \end{array}$$
Suppose $A$,$B$,$C$,...,$a$,$b$,$c$,... are the roots of those two equations. The final equation resulting from the elimination of $y$ must be:
$$\left.\begin{array}{ll} &(A-a)(A-b)(A-c)...\\ \times&(B-a)(B-b)(B-c)...\\ \times&(C-a)(C-b)(C-c)...\\ \times&... \end{array}\right\rbrace=0$$ 
This expression is homogeneous of degree $mn$ in $A$,$B$,$C$,...,$a$,$b$,$c$,..., and symetrical with regard to both sets of roots. We thus could obtain it in terms of the elementary symetrical functions of the roots, \textit{i. e.} in terms of the coefficients $P$,$Q$,...,$p$,$q$,...~; moreover the degree \textit{in $x$} of every coefficient coincides with its degree \textit{in the roots}, so that the degree in $x$ of the final equation is also $mn$.} can only be used on systems of two equations, Bézout adds that Newton's method has the same shortcoming:
\begin{quote}
In fact, Newton's method does not require to compare equations two by two.
Nevertheless, it has no advantage over Euler's and Cramer's method for systems
with more than two equations: then, the final equation is mixed with useless 
factors.\footnote{\textit{cf.} \cite{bezout1764}, p.~290.}
\end{quote}
In 1779, Bézout does not mention explicitely ``Newton's method''. He
rather speaks of ``the principle of substitution''; and the way he uses this 
idea is quite ambiguous. Although the method itself is described as an effective mean of calculating the final equation, it is diverted to the calculation of the \emph{degree} of the final equation. Bézout does never seem to worry about the effective calculation of the final equation.%
\footnote{A simple example could illustrate this problem. Suppose $\deg f^{(2)}=\deg f^{(3)}=2$ and $\deg\phi^{(1)}=4$, and let us make vanish as many terms as possible in $\phi^{(1)}$, thanks to equations $(2)$ and $(3)$. In a naive interpretation of Newton's method, let us start and make vanish terms divisible by $x_3^2$ in $f^{(2)}$ and $\phi^{(1)}$ thanks to $(3)$. If one then makes vanish terms divisible by $x_2^2$ in $\phi^{(1)}$ thanks to $(2)$, other terms divisible by $x_3^2$ will re-appear. A clever use of a ``monomial order'' would remedy this situation, but this was unknown to Bézout. Today, Groebner basis computations rely upon the idea of ordering monomials. About the history of Groebner basis, \textit{cf.} \cite{eisenbud1995} p.~337-338.} After his proof for the case of ``complete equations'', he adds: \label{BezoutPivot}
\begin{quote}
The idea of substitution is the nearest approximation to the elementary ideas of elimination in systems of equations of the first degree\footnote{For $t^{(1)}=t^{(2)}=t^{(3)}=1$, this method embodies what is now termed as ``Gaussian elimination''.}. Although we could apply the same idea to incomplete equations, we are going to present another point of view, that can be applied in a general way, whereas the principle of substitution would need modifications and particular attentions if we should keep up with it.\footnote{\textit{cf.} \cite{bezout1779}, \S~54.}
\end{quote}
Bézout is announcing now a second view on elimination.

\paragraph{Using multiplier-polynomials}
In his second proof\footnote{as he does for ``incomplete equations of the first species'', \textit{cf.} \cite{bezout1779}, \S~59--67.}, Bézout explains that deleting terms in a given polynomial, thanks to equations $(2)$ and $(3)$, amounts to the use of new multiplier-polynomials:
\begin{quote}
  \textit{We ask how many terms we could make vanish in a given polynomial, thanks to these equations, without introducing new terms.}
  
  Suppose that there is only one equation; if, having multiplied it by a polynomial (...), we add the product (...) to the given polynomial: it is obvious that
  \begin{description}
  \item[1\textsuperscript{o}] this addition will not change anything to the value of the given polynomial.
  \item[2\textsuperscript{o}] Supposing the multiplier-polynomial is such as not to introduce new terms, we shall be able to make vanish, in the given polynomial, as many terms as there are in the multiplier-polynomial, because each of them brings one coefficient (...)\footnote{\textit{cf.} \cite{bezout1779}, \S~60.}
  \end{description}
\end{quote}
Thus, in order to make terms vanish in $\phi^{(1)}f^{(1)}$ thanks to $(2)$, Bézout is using a polynomial multiplier $\phi^{(2)}$ of degree $T-t^{(2)}$, and he studies the sum $\phi^{(1)}f^{(1)}+\phi^{(2)}f^{(2)}$.
Then, to make terms vanish thanks to $(3)$, he studies the sum $\phi^{(1)}f^{(1)}+\phi^{(3)}f^{(3)}$. As we can see, the order of proceeding for effective calculation is still imbued with ambiguity. 

\paragraph{The sum-equation}
The first and second views described above are, at best, heuristic views. They could, in no way, be seen as rigorous proofs, nor effective algorithms.
The third view on elimination in Bézout's treatise is explained in book II.
It is the only one that we shall refer to, when summarizing the calculation 
of the degree of the final equation, in next section. At some point before or during the writing of his treatise, Bézout must have become aware of the following fact :  the set of polynomials that are sums of products of $n$ given polynomials by multiplier-polynomials is \emph{closed under this operation}. That is to say, the result obtained after iterating several such operations is again a sum of products of the given polynomials by multiplier-polynomials. All steps become, so to speak, united:
\begin{quote}
From now on, we shall study elimination in a way that differs from what
preceeds, but not essentially.

Let us think that every given equation is multiplied by a specific polynomial,
and that we add up all those products. The result is called ``sum-equation''.
The sum-equation becomes the final equation through the vanishing of all
terms containing any unknown that we should eliminate.

We shall now 1\textsuperscript{o} settle the form of every multiplier-polynomial. 2\textsuperscript{o} Determine how many coefficients, in each multiplier-polynomial, could be considered as useful for the elimination (...)\footnote{\textit{Cf.} \cite{bezout1779}, \S~224.}
\end{quote}
To calculate the degree of the final equation, one thus has to:
\begin{itemize}
\item multiply each equation $f^{(\alpha)}=0$ by a polynomial multiplier $\phi^{(\alpha)}$ with undetermined coefficients, of degree $T-t^{(\alpha)}$;
\item build the ``sum-equation'';
\item ask that all terms vanish, except 1, $x_1,\ x_1^2,...,\ x_1^D$. 
\end{itemize}
In the event of a \emph{single} solution, the method of undetermined coefficients implies that:
$$
\text{nbr. of equations}\geq(\text{nbr. of undetermined coefficients})-(\text{nbr. of useless coeff.})
$$

\section{Bézout's demonstration: second species of incomplete equations}\label{Demonstration2}
We shall now summarize Bézout's demonstration concerning the degree of the eliminand of $n$ incomplete equations of the second species in $n$ unknowns:
$$\left\lbrace\begin{aligned}
f^{(1)}=0\\
f^{(2)}=0\\
\vdots\\
f^{(n)}=0
\end{aligned}\right.$$
where, for all $i$, the polynomial $f^{(i)}$ is a generic member of $C_{\leq t^{(i)},a^{(i)}}$ and the degrees $t^{(i)}$, $a^{(i)}$ verify the restrictive conditions given on page \pageref{secondspecies} above.

To start with, Bézout takes it for granted that there exists a unique ``final equation'' of lowest degree resulting of the elimination of $(n-1)$ unknowns, \textit{i. e.} an eliminand, that could be represented as :
$$\sum_{i=1}^n\phi^{(i)}f^{(i)}=0$$
where the $\phi^{(i)}$ are conveniently chosen ``multiplier-polynomials''
\footnote{\textit{Cf.} \cite{bezout1779} \S 224.}. This being sayed, we are not going to discuss its existence here. The important point is that Bézout is studying the linear map :
$$(\phi^{(1)},\phi^{(2)},...,\phi^{(n)})\longmapsto\sum_{i=1}^n\phi^{(i)}f^{(i)}$$
Doing so, he restricts himself to finite sub-vector spaces by putting an upper bound on the degrees of the $\phi^{(i)}$. For any set of integers $T$, $A_1$, $A_2$,..., $A_n$, $B$ as above, let us call $(f^{(1)},...,f^{(n)})_{\leq T,A,B}$ the linear map defined by
$$\begin{array}{lrcl}
(f^{(1)},...,f^{(n)})_{\leq T,A,B}:&\displaystyle\bigoplus_{i=1}^nC_{\leq T-t^{(i)},A-a^{(i)},B-b^{(i)}}&\longrightarrow &C_{\leq T,A,B}\\
&(\phi^{(1)},\phi^{(2)},..., \phi^{(n)})&\longmapsto &\sum_{i=1}^n\phi^{(i)}f^{(i)}
\end{array}$$
For ease of notation, we shall sometimes write $f_{\leq T,A,B}=(f^{(1)},...,f^{(n)})_{\leq T,A,B}$. We shall also sometimes omit the indices $A$ and $B$ when we fear no confusion. The total number of undetermined coefficients in the $\phi^{(i)}$ polynomials is
$$\dim\bigoplus_{i=1}^nC_{\leq T-t^{(i)},A-a^{(i)},B-b^{(i)}}$$
Bézout says that this number is the number of ``useful coefficients'', plus the number of ``useless coefficients''\footnote{\textit{Cf.} \cite{bezout1779} \S 224.}. In other words :\label{rank}
$$\dim\bigoplus_{i=1}^nC_{\leq T-t^{(i)},A-a^{(i)},B-b^{(i)}}=\dim\text{im}f_{\leq T,A,B}+\dim\ker f_{\leq T,A,B}$$
Now $\text{im} f_{\leq T,A,B}$ is of special interest since any eliminand would belong to it for large enough values of $T$, $A$ and $B$. Moreover, if there exists an eliminand in $x_1$ of lowest degree $D$, then $\lbrace 1,x_1,x_1^2,...,x_1^{D-1}\rbrace$ is a free family in $C_{\leq T,A,B}/\text{im} f_{\leq T,A,B}$. Then one has $D\leq\dim\text{coker} f_{\leq T,A,B}$. We are thus naturally led to calculate\footnote{Let us re-phrase this argument in Bézout's terminology : according to his third view on elimination, the method of undetermined coefficients is leading to the coefficients of the final equation, and these coefficients are the solution of a system of linear equations. One should have, in the event of a single solution :
\begin{align*}
\text{nbr. of undetermined coefficients}&=\dim\bigoplus\limits_{i=1}^nC_{\leq T-t^{(i)},A-a^{(i)},B-b^{(i)}}\\
\text{nbr. of useless coefficients}&=\dim\ker f_{\leq T,A,B}\\
\text{nbr. of linear equations}&=\dim C_{\leq T,A,B}-D\\
\text{nbr. of equations}&\geq(\text{nbr. of undetermined coefficients})-(\text{nbr. of useless coefficients})
\end{align*}
Hence, as written above :
$$D\leq \dim C_{\leq T,A,B}-\dim\bigoplus\limits_{i=1}^nC_{\leq T-t^{(i)},A-a^{(i)},B-b^{(i)}}+\dim\ker f_{\leq T,A,B}$$
As Bézout does never prove the existence of the eliminand, then, what he does actually calculate is the right side of this inequality, \textit{i. e.} $\dim\text{coker} f_{\leq T,A,B}$.
}
\begin{equation}\begin{split}
  \dim\text{coker} f_{\leq T,A,B}=&\dim C_{\leq T,A,B}-\dim\text{im} f_{\leq T,A,B}\\
  =&\dim C_{\leq T,A,B}-\dim\bigoplus\limits_{i=1}^nC_{\leq T-t^{(i)},A-a^{(i)},B-b^{(i)}}+\dim\ker f_{\leq T,A,B}
\end{split}\end{equation}

In \S~233, Bézout describes how to count the number of ``useless coefficients'', \textit{i. e.} $\dim\ker f_{\leq T,A,B}$. He says :

\begin{quote}
  If one remembers what has been said in Book I, one will understand that,
  the number of useful coefficients in the first multiplier-polynomials of
  the equations undergoing elimination, will always be equal to the number
  of coefficients in this polynomial, minus the number of terms that could be
  made to vanish in this polynomial, thanks to the $n-1$ other equations,
  $n$ being the total number of equations;
  
  That the number of useful coefficients in the second multiplier-polynomial,
  will be the total number of coefficients of this polynomial, minus the
  number of terms that could be made to vanish in this polynomial, thanks
  to the $n-2$ last equations;
  
  That the number of coefficients useful in the third multiplier-polynomial,
  will equal the number of terms of this polynomial, minus the number of
  terms that could be made to vanish in this polynomial, thanks to the $n-3$
  other equations; and so on up to the last one, where the number of
  useful coefficients will be precisely equal to the number of
  its terms.\footnote{\textit{Cf.} \cite{bezout1779}, \S~233.}
\end{quote}
The argument is inductive. Let us call $(1),...,(n)$ the $n$ equations.
In order to calculate the number of coefficients that could be made to
vanish in $\phi_1$ using equations $(2),...,(n)$, one should use
new multiplier-polynomials. To paraphrase what Bézout tells us :
\begin{align*}
  \dim(\ker f_{\leq T,A,B}) = & \text{ nbr. of useless coefficients}\\
  = & \text{ nbr. of coeff. to vanish in~}\phi_1\text{ thanks to }(2),...,(n)\\
   & +\text{ nbr. of coeff. to vanish in~}\phi_2\text{ thanks to }(3),...,(n)\\
   & +...\\
   & +\text{ nbr. of coeff. to vanish in }\phi_{n-1}\text{ thanks to }(n)\\
  = & \dim\left(\text{im}(f^{(2)},...,f^{(n)})_{\leq T-t^{(1)},A-a^{(1)},B-b^{(1)}}\right)\\
   & +\dim\left(\text{im}(f^{(3)},...,f^{(n)})_{\leq T-t^{(2)},A-a^{(2)},B-b^{(2)}}\right)\\
   & +...\\
   & +\dim\left(\text{im}(f^{(n)})_{\leq T-t^{(n-1)},A-a^{(n-1)},B-b^{(n-1)}}\right)
\end{align*}
Alas, proving it lies beyond Bézout's means. He seems to have been
aware of the difficulty. The problem could be reduced to proving the following
statement:

\paragraph{Statement}\label{statement} For all $2\leq r\leq n$, for all $(\phi^{(1)},...,\phi^{(r)})$ with 
$$\phi^{(1)}\in C_{\leq T-t^{(1)},A-a^{(1)},B-b^{(1)}},...,\phi^{(r)}\in C_{\leq T-t^{(r)},A-a^{(r)},B-b^{(r)}},$$ such that
$$\sum_{i=1}^r\phi^{(i)}f^{(i)}=0,$$
we have
$$\phi^{(1)}\in\text{im}(f^{(2)},...,f^{(r)})_{\leq T-t^{(1)},A-a^{(1)},B-b^{(1)}}.$$

Although Bézout goes to great lengths studying this situation\footnote{\label{fictitious}See \cite{bezout1779}, \S~107-118, where he tries to convince his reader that it is impossible to increase the number of terms vanishing in $\phi^{(1)}$ by ``fictitious introduction'' (\textit{introduction fictive}) of terms of higher degree.}, there is, as far as we could understand, no proof of this statement in his treatise. For the time being, let us admit this statement (see section \ref{Toric} for the proof). Then we can write:
\begin{align*}
  \dim\ker(f^{(1)},...,f^{(n)})_{\leq T,A,B}=&\dim\text{im}(f^{(2)},...,f^{(n)})_{\leq T-t^{(1)},A-a^{(1)},B-b^{(1)}}\\
  &+\dim\ker(f^{(2)},...,f^{(n)})_{\leq T,A,B}
\end{align*}
By recurrence on the number of equations, we thus have, as written above:
\begin{equation}\begin{split}
  \dim(\ker(f^{(1)},...,f^{(n)})_{\leq T,A,B}) = & \dim\left(\text{im}(f^{(2)},...,f^{(n)})_{\leq T-t^{(1)},A-a^{(1)},B-b^{(1)}}\right)\\
   & +\dim\left(\text{im}(f^{(3)},...,f^{(n)})_{\leq T-t^{(2)},A-a^{(2)},B-b^{(2)}}\right)\\
   & +...\\
   & +\dim\left(\text{im}(f^{(n)})_{\leq T-t^{(n-1)},A-a^{(n-1)},B-b^{(n-1)}}\right)
\end{split}\end{equation}

Let us now use the following notation for finite differences\footnote{
  Bézout's own notation for higher order finite differences could be defined
  by recurrence as follows:
  $$
  d^{r}P(t)\cdots\left(\begin{array}{c} t\\
    -t_{1},...,-t_{r}\end{array}\right)=
  d^{r-1}P(t)\cdots\left(\begin{array}{c} t\\
    -t_{2},...,-t_{r}\end{array}\right)-
  d^{r-1}P(t-t_{1})\cdots\left(\begin{array}{c}t\\
    -t_{2},...,-t_{r}\end{array}\right)
  $$
  The relation between his notation and ours could thus be expressed by:
  $$
  d^{r}P(t)\cdots\left(\begin{array}{c}
t\\
-t_{1},...,-t_{r}\end{array}\right)=\Delta_{t_{1}}...\Delta_{t_{r}}P(t)
  $$
} of a given polynomial $P(T,A,B)$:
$$\Delta_{t,a,b}P(T,A,B)=P(T,A,B)-P(T-t,A-a,B-b)$$

From $(6)$ and $(7)$, by gathering terms, one finds:
\begin{align*}
  \dim\text{coker} f_{\leq T} = & \dim C_{\leq T}-\sum\limits _{i=2}^{n}\dim C_{\leq T-t^{(i)}}+\dim\left(\text{im}(f^{(3)},...,f^{(n)})_{\leq T-t^{(2)}}\right)\\
  &+...+\dim\left(\text{im}(f^{(n)})_{\leq T-t^{(n-1)}}\right)\\
  & -\left(\dim C_{\leq T-t^{(1)}}-\dim\left(\text{im}(f^{(2)},...,f^{(n)})_{\leq T-t^{(1)}}\right)\right)\\
  = & \dim\left(\text{coker}(f^{(2)},...,f^{(n)})_{\leq T}\right)-\dim\left(\text{coker}(f^{(2)},...,f^{(n)})_{\leq T-t_{1}}\right)\\
  = & \Delta_{t^{(1)}}\dim\left(\text{coker}(f^{(2)},...,f^{(n)})_{\leq T}\right)\\
  = & \Delta_{t^{(1)}}\Delta_{t^{(2)}}\dim\left(\text{coker}(f^{(3)},...,f^{(n)})_{\leq T}\right)\\
  = & ...\\
   = & \Delta_{t^{(1)}}\Delta_{t^{(2)}}...\Delta_{t^{(n)}}\dim C_{\leq T}
\end{align*}
where we omit the indices $A,B$ for brevity. For example, the third equality should read:
$$\dim\left(\text{coker}(f_{1},...,f_{n})_{\leq T,A,B}\right)=\Delta_{t^{(1)},a^{(1)},b^{(1)}}\dim\left(\text{coker}(f_{2},...,f_{n})_{\leq T,A,B}\right)$$
This recurrence formula is the heart of Bézout's computations. As $\dim C_{\leq T,A,B}$ is a polynomial of degree $n$ in $T,A,B$, after applying $n$ times the operator $\Delta$, one must obtain a constant independant of $T,A,B$. Eventually, Bézout finds\footnote{More details of this calculation will be given in the proof of prop.~8, section~\ref{Toric}.}:
\begin{multline*}
  \dim\text{coker}(f_{1},...,f_{n})_{\leq T,A,B}\\
  =\prod_{i=1}^nt^{(i)}-\sum_{j=1}^n\prod_{i=1}^n(t^{(i)}-a_j^{(i)})+\prod_{i=1}^n(t^{(i)}-b^{(i)})-\sum_{i=1}^n\left\lbrack(a_1^{(i)}+a_2^{(i)}-b^{(i)})\prod_{j\neq i}(t^{(j)}-b^{(j)})\right\rbrack
\end{multline*}
In 1782, three years after Bézout, Waring takes over this formula in the preface
to the second edition of his \textit{Meditationes algebraicae}\footnote{Waring
  writes (p.~\textit{xvii}-\textit{xx} of the preface to the second edition):
  \begin{quote}
    si sint $(h)$ aequationes $(n,m,l,k,\text{\&c.})$ dimensionum respective
    totidem incognitas quantitates $(x,y,z,v,\text{\&c.})$ involventes~; et
    sint $p$, $q$, $r$, $s$, $\text{\&c.}$~; $p'$, $q'$, $r'$, $\text{\&c.}$~;
    $p''$, $q''$, \&c.~; maximae dimensiones, ad quas ascendunt
    incognitae quantitates $x$, $y$, $z$, $v$, \&c., in respectivis
    aequationibus $(n,m,l,k,\text{\&c.})$ dimensionum~; tum aequationem, cujus
    radix est $x$ vel $y$ vel $z$, \&c. haud ascendere ad plures quam
    $n\times m\times l\times k\times\text{\&c.}-(n-p)\times(m-p')\times(l-p'')\times\text{\&c.}-(n-q)\times(m-q')\times(l-q'')\times\text{\&c.}-(n-r)\times(m-r')\times\text{\&c.}-\text{\&c.}=P$
    dimensiones~: si vero dimensiones quantitatum $(x \text{ \& } y)$ simul
    sumptarum haud superent dimensiones $a,a',a'',\text{\&c.}$ in predictis
    aequationibus~; tum aequationem, cujus radix est $x$, \&c. haud plures
    habere quam $P+(n-a)\times(m-a')\times\text{\&c.}-(p+q-a)\times(p'-a')\times\text{\&c.}$
  \end{quote}
  Waring forgets to mention the ``restrictive conditions'' (see
  p.~\pageref{secondspecies} above). With his notations, those conditions require that:
  $$\begin{array}{llll}
    r+a\geq n & r'+a'\geq m & r''+a''\geq l & ...\\
    s+a\geq n & s'+a'\geq m & s''+a''\geq l\\
    \vdots & \vdots & \vdots
  \end{array}$$
}. 
Most importantly, in 1779, Bézout had also noticed that this $n$-th order finite difference could be written as an alternate sum\footnote{\textit{Cf.} \cite{bezout1779} p.~43.}:\label{alternatesum}
$$
\Delta_{t^{(1)}}\Delta_{t^{(2)}}...\Delta_{t^{(n)}}\dim C_{\leq T}=\dim C_{\leq T}-\dim\bigoplus_{i}C_{\leq t-t^{(i)}}+\dim\bigoplus_{i<j}C_{\leq T-t^{(i)}-t^{(j)}}-...
$$
Later, Cayley will shed light upon this alternate sum.

\section{Complete equations and the first species of incomplete equations}

As was said above, one could, from the formula for the degree of the eliminand of $n$ incomplete equations of the second species, derive the formulae for complete equations and for the first species of incomplete equations.
For $n$ incomplete equations of the first species, \textit{i. e.} equations in $x_1,x_2,...,x_n$ of the form:
$$
\sum\limits_{\begin{array}{l}
\scriptstyle k_{1}\leq a_{1},\ k_{2}\leq a_{2},\ k_{3}\leq a_{3},...\\
\scriptstyle k_{1}+k_{2}+k_{3}+...\leq t\end{array}}{u_kx^k}=0
$$
where
$$\left\lbrace\begin{array}{l}
(\forall i)\ a_i\leq t\\
(\forall i\neq j)\ a_i+a_j>t
\end{array}\right.$$
the degree of the final equation resulting from the elimination of $x_{2},x_{3},...$ is:\footnote{This formula is still quoted in 1900 in Netto's \textit{Vorlesungen über Algebra} \cite{netto1900}, \S~419.}
$$D=\prod_{i=1}^nt^{(i)}-\sum_{i=1}^n\prod_{j=1}^n(t^{(j)}-a_i^{(j)})$$
Now specifying, for all $i,j$, $a_j^{(i)}=t^{(i)}$, one obtains the well-known formula for the case of $n$ ``complete equations'', \textit{i. e.} generic of their degrees :
$$D=\prod_{i=1}^nt^{(i)}$$
In fact, for general $n$, only this case of Bézout's theorem will be the
object of rigorous study in XIXth century (\textit{cf.} Serret \cite{serret1866}, Schmidt \cite{schmidt1886}, Netto \cite{netto1900} vol.~2).

\section{Bézout's theorem for three incomplete equations of the third species}\label{Demonstration3}

As for the ``third species of incomplete equations'', when $n=3$, Bézout is hitting another major problem : $\dim C_{\leq T,A,B}$ is not any more a polynomial in $T,A,B$. It is, so to speak, polynomial by pieces. Let us write
$$H_1=T-B_2-B_3+A_1,\qquad H_2=T-B_3-B_1+A_2,\qquad H_3=T-B_1-B_2+A_3.$$
For each of the eight ``forms'' corresponding to different algebraic signs of those three quantities, $\dim C_{\leq T,A,B}$ is a different polynomial in $T,A,B$. More precisely, calling $P_i(T,A,B)=\dim C_{\leq T,A,B}$ when the values of $T,A,B$ belong to the $i$-th form, one has:
\begin{description}
\item[1st form] ($H_1\leq0,\quad H_2\leq0,\quad H_3\leq0$):
  \begin{align*}
    P_1(T,A,B) = & {3+T \choose 3}-\sum_{i=1}^3{3+T-A_i-1 \choose 3}+\sum_{i=1}^3{3+T-B_i-2 \choose 3}\\
    & -\sum_{i=1}^3\left[(A_1+A_2+A_3-A_i-B_i){2+T-B_i-1 \choose 2}\right]
  \end{align*}
\item[2d form] ($H_1\leq0,\quad H_2\leq0,\quad H_3\geq0$):
  $$P_2(T,A,B)=P_{1}(T,A,B)+{3+T+A_3-B_1-B_2-2 \choose 3}$$
\item[3rd form] ($H_1\geq0,\quad H_2\leq0,\quad H_3\leq0$):
  $$P_3(T,A,B)=P_{1}(T,A,B)+{3+T+A_1-B_2-B_3-2 \choose 3}$$
\item[4th form] ($H_1\geq0,\quad H_2\leq0,\quad H_3\geq0$):
  $$P_4(T,A,B)=P_{3}(T,A,B)+{3+T+A_3-B_1-B_2-2 \choose 3}$$
\item[5th form] ($H_1\geq0,\quad H_2\geq0,\quad H_3\leq0$):
  $$P_5(T,A,B)=P_{3}(T,A,B)+{3+T+A_2-B_1-B_3-2 \choose 3}$$
\item[6th form] ($H_1\geq0,\quad H_2\geq0,\quad H_3\geq0$):
  $$P_6(T,A,B)=P_{5}(T,A,B)+{3+T+A_3-B_1-B_2-2 \choose 3}$$
\item[7th form] ($H_1\leq0,\quad H_2\geq0,\quad H_3\leq0$):
  $$P_7(T,A,B)=P_{1}(T,A,B)+{3+T+A_2-B_1-B_3-2 \choose 3}$$
\item[8th form] ($H_1\leq0,\quad H_2\geq0,\quad H_3\geq0$):
  $$P_8(T,A,B)=P_{7}(T,A,B)+{3+T+A_3-B_1-B_2-2 \choose 3}$$
\end{description}
Suppose that the argument developed above for the second species of incomplete equations is also valid for the third species, then one should have, as above:
$$\dim\text{coker}(f^{(1)},...,f^{(n)})_{\leq T,A,B}=\Delta_{t^{(1)},a^{(1)},b^{(1)}}\Delta_{t^{(2)},a^{(2)},b^{(2)}}...\Delta_{t^{(n)},a^{(n)},b^{(n)}}\dim C_{\leq T,A,B}$$
The rest of the computation could only be done under the assumption that all vector-spaces $C_{\leq...}$ actually involved in this expression belong to the same ``form''. Let us write
$$D_i=\Delta_{t^{(1)},a^{(1)},b^{(1)}}\Delta_{t^{(2)},a^{(2)},b^{(2)}}...\Delta_{t^{(n)},a^{(n)},b^{(n)}}P_i(T,A,B)$$
Bézout calculates $D_1,D_2,...,D_8$. The actual eight formulae occupy no less than eight full pages of the treatise\footnote{\textit{Cf.} \cite{bezout1779} \S 119-127.}. He then proposes a test, or rather, ``symptoms''\footnote{See the ``symptoms enabling us to recognize, among the different expressions of the value of the degree of the final equation, those that one should choose or reject'', in \cite{bezout1779} \S 117. Here again, Bézout's own justifications lack evidence; they rely upon undemonstrated facts about the sum-equation, as in footnote \ref{fictitious} p.~\pageref{fictitious} above.} to reject some of those eight values. The other values are ``admissible''; as such, all of them must be equal to the degree of the eliminand, according to Bézout. In section \ref{Toric} below, we shall prove that Bézout's choice was right.

\section{The theory of the resultant in XIXth century}

We have just presented Bézout's slowly maturated treatise and the fertile historical context of its publication. We are struck by the lack of immediate posterity of this book: sixty years separate the publication of Bézout's treatise and the first revival of what was reckoned, in XIXth century, as the ``theory of elimination''. Bézout's treatise was complex, it was clearly perceived as such by one of his early readers, and this fate follows it until today. Poisson recognized the importance of Bézout's work but immediately pointed out to the gap between the strength of the theorem and the ``difficulties'' of its demonstration :
\begin{quote}
This important theorem is Bézout's, but the way he proves it is neither direct nor simple; nor is it devoid of any difficulty.\footnote{\textit{Cf.} \cite{poisson1802} p. 199. See also Brill and Noether, saying about Bézout's book that it is ``as well-known as lacking readers'', and that, by the time of Jacobi, ``most of it had fallen into oblivion'', \textit{cf.} \cite{brill1892} p. 143 and 147.}
\end{quote}
Three mathematiciens produced, so to speak, a new beginning in elimination theory, between 1839 and 1848: their names are Sylvester, Hesse, Cayley\footnote{Some historians and mathematicians have said that Sylvester and Richelot (1808-1875, student of Jacobi) had discovered the ``dialytic'' method of elimination, although this method was not very different from Bézout's method when dealing with only two equations ; it is also said that Hesse had re-discovered this method, in 1843. \textit{Cf.} Max Noether, \cite{noether1898}, p.~136, about Sylvester \cite{sylvester1839}, \cite{sylvester1840} and \cite{sylvester1841}, and Richelot \cite{richelot1840}. Eberhard Knobloch \cite{knobloch1974,knobloch2001} noticed that Leibniz already knew this method ; perhaps Leibniz was even closer to Sylvester's ideas, than Euler and Bézout. In what follows, we shall only draw a comparison between the works of Sylvester, Richelot and Hesse, and those of Bézout from 1779 concerning the \emph{sum-equation}, with an emphasis on the case of $n$ equations when $n>2$.}. The main object of study is, rather than the eliminand, the ``resultant'' of $n$ homogeneous polynomials in $n$ indeterminates.

Before entering into the works of those three scholars, it is to be noted that two special cases progressed in first half of XIXth century: linear equations, thanks to the theory of determinants, and the case of two equations in one or two unknowns. When eliminating one unknown between two equations in two unknowns, one obtains the eliminand. When eliminating one unknown between two equations in one unknown (or between two \emph{homogeneous} equations in two unknowns, which is the same thing), one obtains the resultant. The discriminant is the resultant of a polynomial and its derivate. The interest in the discriminant is motivated by the study of the euclidian algorithm for polynomials, following Sturm's researches about the roots of polynomials over $\mathbb{R}$. The case of two equations also benefits from the methods of analysis, for example in Ossian Bonnet's works culminating in 1847 when he defines the intersection multiplicity of two curves in one point.

\paragraph{Sylvester} Sylvester's researches are stimulated by Sturm's theorem. When applying the euclidian division algorithm to two polynomials in one indeterminate, the successive remainders are also called ``Sturm functions''. In 1840, Sylvester gives a formula in terms of determinants to calculate Sturm functions. As was known long before, the last remainder, being of degree 0, can be seen as the result of the elimination of the indeterminate. This is what Sylvester calls the ``dialytic method of elimination''. There is no demonstration in this short article. Maybe Sylvester did not know, in 1840, of Euler's and Bézout's works about elimination. He does not refer to them ; later, in 1877, he himself says that he had discovered the dialytic method by teaching to a pupil\footnote{In \cite{sylvester1877}, Sylvester says: ``I remember,...,
how... when a very young professor, fresh from the University of Cambridge,
in the act of teaching a private pupil the simpler parts of Algebra,
I discovered the principle now generally adopted into the higher text
books, which goes by the name of the \textit{Dialytic Method of Elimination}''.}. 

In 1841, Sylvester develops the dialytic method in the case of three quadratic \emph{homogeneous} equations in three unknowns $x$, $y$, $z$. The interest in homogeneous equations is crucial to our subject, and it could well be explained by the fusion between projective geometry and algebraic geometry under the influence of Möbius and Plücker. Sylvester delelops several versions of his dialytic method. In one of them\footnote{\textit{Cf.} \cite{sylvester1841}, example 4, p.~64; other versions are given in footnotes.}, he mutiplies each equation by the three monomials of degree 1. Thus, for the system
$$\left\lbrace\begin{array}{l}
U=0\\
V=0\\
W=0
\end{array}\right.$$
one obtains 9 equations of degree 3 :
$$
\left\lbrace \begin{array}{lll}
xU=0 & yU=0 & zU=0\\
xV=0 & yV=0 & zV=0\\
xW=0 & yW=0 & zW=0\end{array}\right.
$$
As there exists 10 monomials of degree 3, we are short of one equation to apply the dialytic method and build a determinant. Sylvester uses the jacobian determinant to get a tenth equation of degree 3:
$$
\frac{1}{8}\times\begin{vmatrix}\dfrac{\partial U}{\partial x} & \dfrac{\partial U}{\partial y} & \dfrac{\partial U}{\partial z}\\[1em]
\dfrac{\partial V}{\partial x} & \dfrac{\partial V}{\partial y} & \dfrac{\partial V}{\partial z}\\[1em]
\dfrac{\partial W}{\partial x} & \dfrac{\partial W}{\partial y} & \dfrac{\partial W}{\partial z}\end{vmatrix}=0
$$
One thus obtains 10 equations in the 10 monomials of degree 3. If one considers those monomials as the 10 independant unknowns of a system of 10 linear equations, elimination is reduced to the calculation of a single determinant. As compared to the method of the sum-equation, this method is sparing the use of undetermined coefficients. Moreover, it is a \emph{symbolic} method. It uses the ambivalence of the symbolic expression of the monomials: each monomial is both a monomial in the unknowns of the initial system, and an independant unknown of a system of linear equations. This allows the transfer of the symbolic methods of linear algebra (determinants, and soon, matrices) to the algebra of forms of higher degree\footnote{Sylvester is probably refering to this in particular when he says that the ``great principle of dialysis, originally discovered in the theory of elimination, in one shape or another pervades the whole theory of concomitance and invariants'', \textit{cf.}  \cite{sylvester1852a}, p.~294.}.

Contemporary readers must have been surprised by the use of the jacobian determinant. In general, for equations of higher degree and systems with more unknowns, there still remains to explain the appropriate choice of linear equations. In an article \cite{sylvester1841b} published in the same year, Sylvester gives a general method. We are not going to describe it entirely. Suffice to say that, after having obtained a first set of equations called \emph{augmentatives} by multiplying each initial equation by the monomials (like the 9 equations of degree 3 above), Sylvester builds other equations called \emph{secondary derivatives}, as follows. He writes each of the $n$ initial equations under the form $x^{\alpha}F+y^{\beta}G+z^{\gamma}H+...$ where $x$, $y$, $z$,... are the $n$ unknowns to eliminate, $F$, $G$, $H$,... are polynomials, and $\alpha$, $\beta$, $\gamma$,... is any system of integers allowing such a representation. He thus obtains a system of $n$ equations:
$$
\left\lbrace \begin{array}{l}
x^{\alpha}F+y^{\beta}G+z^{\gamma}H+...=0\\
x^{\alpha}F'+y^{\beta}G'+z^{\gamma}H'+...=0\\
\vdots\end{array}\right.
$$
The following determinant is one of the secondary derivatives:
$$
\begin{vmatrix}
F & G & H & \cdots\\
F' & G' & H' & \cdots\\
\vdots\end{vmatrix}=0
$$
The many choices possible for $\alpha$, $\beta$, $\gamma$, \textit{etc.} allow as many secondary derivatives. 

When $n=3$ et $\deg U=\deg V=\deg W=m$, taking $\alpha=\beta=\gamma=1$ and
\begin{align*}
F=\frac{1}{m}\frac{\partial U}{\partial x},\quad G=\frac{1}{m}\frac{\partial U}{\partial y},\quad H=\frac{1}{m}\frac{\partial U}{\partial z}\\
F'=\frac{1}{m}\frac{\partial V}{\partial x},\quad G'=\frac{1}{m}\frac{\partial V}{\partial y},\quad H'=\frac{1}{m}\frac{\partial V}{\partial z}\\
F''=\frac{1}{m}\frac{\partial W}{\partial x},\quad G''=\frac{1}{m}\frac{\partial W}{\partial y},\quad H''=\frac{1}{m}\frac{\partial W}{\partial z}
\end{align*}
one finds the jacobian determinant, up to a constant factor.

Let us now observe the easy case $n=2$, $m=\deg U=\deg V$. Sylvester does not even mention it in his article. If $1\leq\alpha\leq m$ and $\beta=0$, taking as $G$ (resp. $G'$) the sum of terms of $U$ (resp. $V$) of degree less than $\alpha$ in $x$, one has, for each value of $\alpha$, one equation of degree $(m-1)$ in $x$ :
$$
\begin{vmatrix}
F & G\\
F' & G'\end{vmatrix}=0
$$
Hence there is no need of augmentatives. The expression obtained by eliminating dialytically between the $m$ equations of degree $(m-1)$ is none other than the determinant of a matrix that Bézout had already studied \cite{bezout1764} in 1764. Bézout's matrix might thus have inspired this general method to Sylvester. This matrix also played an important role in his famous later memoir \textit{On a theory of syzygetic relations}\footnote{\textit{Cf.} \cite{sylvester1853a}.}.

Still in 1841, Sylvester also considers the possibility of building augmentatives of degree at least $\sum t_{i}-n+1$, where $t_{i}$ is the degree of the $i$-th initial equation. In this case, the augmentatives suffice to build a determinant, and there is no need of secondary derivative. This is the grounding for a method developed by Cayley in 1848 (see below).

One must say that Sylvester's elimination method often brings out a superfluous factor; Sylvester doesn't give any mean of detecting and isolating this factor. His method leads directly to the resultant in but a few cases, as the two cases mentioned above for two or three equations. Despite of this drawback, Sylvester has clearly circumscribed a domain of research not limited to the resultant or the eliminand.\footnote{In 1997, Jean-Pierre Jouanolou gave a complete study of the secondary derivatives that he calls ``formes de Sylvester'', in \cite{jouanolou1997}, \S~3.10.}

\paragraph{Hesse}
Applying elimination to the study of plane cubics in 1844, Hesse goes back to the formalism of Bézout's ``sum-equation''. He knew of Bézout's treatise, and he does mention it\footnote{\textit{Cf.} \cite{hesse1844}. Hesse also knew of Richelot and Sylvester, and of the works of Euler. He extends to $n>2$ equations the method of Euler-Cramer using symetrical functions, as Poisson had done in 1802, although he probably didn't know of Poisson's article.}. Let there be three quadratic equations in two unknowns:
$$\left\lbrace\begin{array}{l}
U=0\\
V=0\\
W=0
\end{array}\right.$$
In order to eliminate the two unknowns, Bézout would have used multiplier-polynomials of degree 2. Following the way of thinking of XIXth century algebraists, let us homogenize the sum-equation of degree 4, thus getting a third unknown $z$. Then there exist multiplier-polynomials $A$, $B$, $C$ such that:
$$AU+BV+CW=z^4R$$
where $R=0$ is the resultant of $U$, $V$, $W$\footnote{It is to be noted that Sylvester had also mentioned this sum-equation, translated in his own dialytic formalism where one would multiply by monomials of degree 2, in a footnote in \cite{sylvester1841}, p.~64-65.}. We must keep in mind this sum-equation when studying the other equations derived by Hesse:
\begin{itemize}
\item First of all, Hesse translates as a sum-equation the method of Sylvester using the jacobian determinant. If $\phi$ is the jacobian determinant of $U$, $V$, $W$, one could obtain a sum-equation of degree 3 thanks to multiplier-polynomials of degree 1:
$$AU+BV+CW+\delta\phi=z^3R$$
    where $\deg A=\deg B=\deg C=1$ and $\delta$ is a constant\footnote{\textit{Cf.} \cite{hesse1844}, \S~8-10.}.
\item Hesse then observes that the jacobian itself can be obtained by a sum-equation with multiplier-polynomials of degree 2:
      $$AU+BV+CW=z\phi$$
\item He also proves that the partial derivatives of $\phi$ could be obtained in the same fashion, under the form
$$AU+BV+CW+\delta\phi=z\frac{\partial\phi}{\partial x}$$
\item One is thus allowed to calculate the resultant $R$ with multiplier-polynomials of degree 0, thanks to the partial derivatives of the jacobian determinant
\footnote{\textit{Cf.} \cite{hesse1844}, \S~11-14.}:
$$aU+bV+cW+d\frac{\partial\phi}{\partial x}+e\frac{\partial\phi}{\partial y}+f\frac{\partial\phi}{\partial z}=z^2R$$
\end{itemize}
As modern elimination-theorists would say\footnote{As we shall see further, with Hurwitz \cite{hurwitz1913}, one defines the ideal of \emph{Trägheitsformen} of $(U,V,W)$ as
  $$\mathfrak{m}^{-\infty}(U,V,W)=\bigcup\limits_{k\geq0}\left\lbrace f\ \mid\ \mathfrak{m}^kf\subset(U,V,W)\right\rbrace$$
where $\mathfrak{m}=(x,y,z)$. This langage is still in use today, \textit{cf.} \cite{jouanolou1991}.}, Hesse thus studied the resultant, the jacobian, and the partial derivatives of the jacobian, as \emph{Trägheitsformen}, or inertia forms, of the ideal $(U,V,W)$. We won't describe how Hesse applied these calculations to the case where $U$, $V$, $W$ are the partial derivatives of the homogeneous polynomial defining a plane cubic.

\paragraph{Cayley} The concept of ``sum-equation'' appears again in 1847, in an article by Cayley, \textit{On the theory of involution in geometry}. Cayley says that a homogeneous polynomial $\Theta$ of degree $r$ is ``in involution'' with homogeneous polynomials $U$, $V$, ... of degrees $m,n,...$ if
$$\Theta=AU+BV+...,$$
\textit{i. e.} if $\Theta\in(U,V,...)_r$, where $(U,V,...)_r$ is the componant of degree $r$ in the homogeneous ideal $(U,V,...)$. He says that ``there is also an analytical application of the theory, of considerable interest, to the problem of elimination between any number of [homogeneous] equations containing the same number of variables''. He conceives of this elimination as the result of Sylvester's dialytic method, but he traces back the origin of his research to ``Cramer's paradox'' and related works by Euler, Cramer, Plücker, Jacobi and Hesse. In this article and another one of 1848, according to some mathematicians of our day, ``Cayley in fact laid out the foundations of modern homological algebra''\footnote{\textit{Cf.} \cite{gelfand1994} p.~\textit{ix}.}. 

Cayley asks for the number of ``arbitrary constants'' in $\Theta$. In modern language, this number is
$$\dim(U,V,...)_r$$
Let us follow Cayley's application of the dialytic method. Let us multiply $U$ by each monomial of degree $r-m$ (for large enough $r$), and $V$ by each monomial of degree $r-n$, and so on. Let $C$ be the ring of polynomials. We shall use the compact notation :
$$\sum\dim C_{r-m}=\dim C_{r-m}+\dim C_{r-n}+...$$
This number is thus the number of linear equations to be solved, when using the dialytic method. The $\dim(C_r)$ monomials of degree $r$ will be the independant unknowns of this system of linear equations. We shall represent this system by a matrix where each column corresponds to one equation\footnote{The matrices mentioned here and below, and the solution of this problem of linear algebra, only appeared in a second article \cite{cayley1848} published by Cayley in 1848. Moreover, we should rather speak of ``matricial figure'' than ``matrix'', because this mathematical object was not yet seen as an \emph{operator}, and its properties were not yet completely developed in the 1840's. The figures 1, 2 and 3, at the end of our article, are lose reproduction of figures in \cite{cayley1848}.}. This matrix (\textit{cf.} figure 1 at the end of this article) is the matrix of a linear map :
$$f_1:\bigoplus C_{r-m}\longrightarrow C_r$$
Going back to the concept of ``sum-equation'' or ``involution'' as Cayley used to say, one sees that $f_1$ is defined by :
$$f_1(A,B,...)=AU+BV+...$$

Cayley notes that the columns of the matrix are not linearly independant, and he asks how many independant columns there are, \textit{i. e.}
$$\dim(\text{im}f_1)$$
In order to eliminate dialytically, one needs $\dim C_r$ independant columns ; then, one could extract a square sub-matrix and build the determinant. To check that the elimination is possible, Cayley sets out to calculate
$$N=\dim C_r-\dim(\text{im}f_1)$$
The relations of linear dependancy between columns are given by families of coefficients that Cayley writes as \emph{lines} of coefficients \emph{under} the original matrix (\textit{cf.} figure 2). In terms of the sum-equation, these relations constitute the kernel of $f_1$. Cayley admits \emph{without proof}\footnote{\textit{Cf.} \cite{cayley1847a}, p.~261, ``the number $N$ must be diminished
by...''.} that $\ker f_1$ is generated by elements of $\bigoplus C_{r-m}$ of the following form:
$$
\left\lbrace \begin{array}{ll}
(MV,-MU,0,...,0)\quad & \text{where}\quad M\in C_{r-m-n}\\
(MW,0,-MU,0,...,0)\quad & \text{where}\quad M\in C_{r-m-p}\\
\vdots\end{array}\right\rbrace
$$
Hence, it is generated by $\sum\dim C_{r-m-n}$ vectors. The height of the second bloc in the matricial figure is thus $\sum\dim C_{r-m-n}$, and $\ker f_1$ is generated by the image of a second map $f_2$ :
$$\begin{array}{llll}
f_{2}: & \bigoplus C_{r-m-n} & \longrightarrow & \bigoplus C_{r-m}\\
 & (M,0,...,0) & \longmapsto & (MV,-MU,0,...,0)\\
 & (0,M,0,...,0) & \longmapsto & (MW,0,-MU,0,...,0)\\
& ...
\end{array}$$

By iterating this procedure, Cayley obtains a figure that we reproduce on figure 3 ; but he does not explain the following steps. Following the path suggested by Cayley, one should endeavour to write the sequence of linear maps thus obtained, and prove that it is an exact sequence\footnote{\textit{Cf.} \cite{cayley1848}, p.~371. Of course, the fact that the sequence is exact is not at all obvious, and Cayley did not prove it. The first historical proof will be mentioned in section \ref{koszul} below.}:
$$...\longrightarrow\bigoplus C_{r-m-n-p}\xrightarrow{f_{3}}\bigoplus C_{r-m-n}\xrightarrow{f_{2}}\bigoplus C_{r-m}\xrightarrow{f_{1}}C_{r}$$
where $\text{im}f_{i+1}=\ker f_i$. Cayley concludes rashly:
\begin{align*}
N & = \dim C_r-\dim(\text{im} f_1)\\
 & = \dim C_r-\sum\dim C_{r-m} +\sum\dim C_{r-m-n}-\sum\dim C_{r-m-n-p} +...
\end{align*}
For large enough values of $r$, this quantity can be expressed in terms of binomial coefficients. As a matter of fact, $N=0$ : hence, according to Cayley, dialytic elimination is possible, for any given number of homogeneous equations with as many unknowns.

What also distinguishes XIXth century authors from their predecessors including Bézout, is the endeavour to give explicit formulae for the result of elimination. Cayley is looking for a formula of the resultant. In 1847, he gives the following formula:
$$R=\frac{Q_{1}Q_{3}...}{Q_{2}Q_{4}...}$$
where, for all $i$, $Q_i$ is a subdeterminant from the matrix of $f_i$. The choice of these subdeterminants obeys the following rule. There are as many columns in the matrix of $f_i$, as lines in the matrix of $f_{i+1}$. The rule is that the set of columns occurring in $Q_i$ must be the complement of the set of lines occurring in $Q_{i+1}$. As for the last one, say $Q_j$, it is the determinant of a maximal square submatrix of the matrix of $f_j$, and it must be chosen such that all $Q_i$ are non-zero. Several authors (Salmon, Netto) seem to have tried proving Cayley's formula, until Macaulay resigned this task and found another formula, simple but ``less general'' \label{FormuleQuotient}, of the form $\dfrac{Q_1}{\Delta}$ where $Q_1$ is still a subdeterminant of the matrix of $f_1$ although its choice is more constrained\footnote{\textit{Cf.} \cite{macaulay1903}, and also \cite{jouanolou1997} for a recent study of Macaulay's formulae.}. En 1926, E. Fischer achieved demonstrating Cayley's formula\footnote{\textit{Cf.} \cite{piel1934}.}.

\section{Koszul complex}\label{koszul}

We have seen the role of an exact sequence leading to an alternate sum of dimensions in Cayley's work about the resultant, although Cayley describes only the two first maps of the sequence, and he omits any proof of exactness.

On the other hand, let us look back at Bézout's work on complete equations. Bézout's alternate sum\footnote{\textit{Cf. supra} p. \pageref{alternatesum}.} obtained through finite differences is quite similar to Cayley's alternate sum. There is just one more unknown in the equations given, because Bézout is studying the eliminand, whereas Cayley is studying the resultant. Bézout's result could be deduced from Cayley's through homogenizing.

As mentioned previously, Serret and Schmidt gave a rigorous proof of Bézout's theorem in the case of the eliminand of $n$ complete equations in $n$ unknowns. Their method could also apply to the case of the resultant. But Serret and Schmidt bypass the need to study the exact sequence above: their proof is an \textit{a posteriori} proof of the degree of the eliminand.

\paragraph{Hurwitz} The first in-depth study of the two first maps of the sequence described by Cayley was given by Hurwitz in 1913. Hurwitz follows the ideas of Mertens \cite{mertens1886} who had already given a complete but complicated theory of the resultant in 1886. In this paragraph and the next, we shall consider $r$ homogeneous polynomials $f_1,f_2,...,f_r$ of degrees $t_1,t_2,...,t_r$ in a polynomial ring $C=K[a,b,...][X_0,...,X_n]$ over a number field, where $X_0,...,X_n$ on one hand, $a,b,...$ on the other hand, are indeterminates. The indeterminates $a,b,...$ will serve the purpose of building ``generic'' polynomials in $X_0,...,X_n$. Let us call $a_{\alpha k}$ the coefficient of $X_k^{t_\alpha}$ in $f_\alpha$. Consider the ideal $\mathfrak{M}=(X_0,...,X_n)$. Hurwitz introduces a new concept : he says that $f\in C$ is a \textit{Trägheitsform} if there exists un integer $k\geq 0$ such that $\mathfrak{M}^kf\subset(f_1,...,f_r)$. For all $k$, let us write, as Hurwitz does, $[f]_k$ the result of substituting $\dfrac{a_{\alpha k}X_{k}^{t_{\alpha}}-f_{\alpha}}{X_{k}^{t_{\alpha}}}$ to $a_{\alpha k}$ for all $\alpha$ in $f$. In the case of generic homogeneous polynomials $f_1,...,f_r$, \textit{i. e.} when their coefficients are the indeterminates $a,b,...$, the following propositions are equivalent%
\footnote{\textit{Cf.} \cite{hurwitz1913}; Hurwitz proves $(\mathrm{i})\iff(\mathrm{ii})$ in his proposition 1, and $(\mathrm{ii})\iff(\mathrm{iii})\iff(\mathrm{iv})\iff(\mathrm{v})$ in proposition 9. It is then obvious that $(\mathrm{iv})\iff(\mathrm{v})\iff(\mathrm{vi})$. The importance of criterion $(\mathrm{iii})$ appears in Mertens' work. For $(\mathrm{ii})\iff(\mathrm{vi})$, see also Jouanolou \cite{jouanolou1991}, (4.2.3) p.~132.}:
\begin{description}
 \item[$(\mathrm{i})$] for $t$ large enough, $(f_1,...,f_r)_t=(f,f_1,...,f_r)_t$ 
 \item[$(\mathrm{ii})$]  $f$ is a \textit{Trägheitsform}
 \item[$(\mathrm{iii})$]  there exists $k$ such that $[f]_k$ is zero
 \item[$(\mathrm{iv})$]  there exists $k,m$ such that $X_k^mf\in(f_1,...,f_r)$
 \item[$(\mathrm{v})$] for all $k$ there exists $m$ such that $X_k^mf\in(f_1,...,f_r)$
 \item[$(\mathrm{vi})$] there exists $m$ such that $X_0^mf\in(f_1,...,f_r)$
\end{description}
The smallest integer $k$ such that $\mathfrak{M}^kf\subset(f_1,...,f_r)$
is called the ``rank'' of $f$ (\textit{Stufe}). A \textit{Trägheitsform}
is said ``proper'' if its rank is non zero. Hurwitz gives a general study
of the \textit{Trägheitsformen}. Let us first consider the case where $n=r$.
Hurwitz proves that the \textit{Trägheitsformen} of rank $\sigma$ are of degree
$\sum t_{\alpha}-n-\sigma+1$. He shows that all \textit{Trägheitsformen} of
rank 1 belong to the ideal $(J,f_1,...,f_n)$ where $J$ is the jacobian
determinant of $f_1,...,f_n$; in rank $>1$, he succeeds in giving
explicit formulae at least for some of the \textit{Trägheitsformen},
those that are linear in the coefficients of each of
the polynomials $f_{1},...,f_{n}$. He proves the existence of
a \textit{Trägheitsform} of degree 0 which generates
the ideal of all \textit{Trägheitsformen} of degree 0: this is the resultant.
In the case $r<n$, Hurwitz succeeds in proving that $(f_1,...,f_r)$ has no
proper \textit{Trägheitsform}. Consider the following assertions:
\begin{description}
\item[$(\mathrm{I}_{n})$] for all $r<n$, the ideal $(f_{1},...,f_{r})$ has no proper Trägheitsform
\item[$(\mathrm{II}_{n})$] for all $r\leq n$, if $A_{1}f_{1}+...+A_{r}f_{r}=0$, then there exists some $L_{ij}$ such that $(\forall i)\ A_{i}=\sum L_{ij}f_{j}$ and $(\forall i,j)\ L_{ij}=-L_{ji}$
\end{description}
Hurwitz proves $\mathrm{I}_n\Rightarrow \mathrm{II}_n\Rightarrow \mathrm{I}_{n+1}$. We recognize in $\mathrm{II}_n$ the assertion made by Cayley about the exactness of his sequence in degree 1. Finally, without any condition on $r$, $n$, Hurwitz also proves, with a similar argument, that $(f_1,...,f_r)$ has no proper \textit{Trägheitsform} of degree $>\sum t_\alpha-r$.

\paragraph{Koszul's complex} An explicit description of the exact sequence touched upon by Bézout, and later conjectured by Cayley, is to be found in the work of the algebraist Koszul taking over the tools of differential geometry in the late 1940's. Koszul's complex is an avatar of de Rham's complex of differential forms\footnote{Retrospective studies on Koszul's work are to be found in \textit{Annales de l'Institut Fourier}, 37 (1987). See the allocution by H. Cartan \cite{cartan1987}.}. Let there be given $r$ homogeneous polynomials $f_1,f_2,...,f_r\in C=K[X_0,...,X_n]$. Suppose that for all $r\leq n$, $f_r$ does not divide zero modulo $(f_1,...,f_{r-1})$. This hypothesis\footnote{\textit{Cf.} \cite{serre2000}, p. 59, where Serre calls \textit{$M$-sequence} such a sequence of polynomials. Some authors speak of a \textit{regular sequence} (\textit{cf.} \cite{hartshorne1977}, II.8, p.~184). See \cite{bourbaki1958} for the main properties of Koszul's complex.} follows from property $(\mathrm{II}_n)$ above. Under this hypothesis, there exists a free resolution of $C/(f_1,f_2,...,f_n)$ bearing Koszul's name. Put $M=\bigoplus\limits_{i=1}^rCe_i$ the free $C$-module of rank $r$. Koszul's complex is the sequence of maps:
$$0\longrightarrow C=\Lambda^0M\longrightarrow\Lambda^1M\longrightarrow\Lambda^2M\longrightarrow\cdots\longrightarrow\Lambda^{r-1}M\longrightarrow\Lambda^rM\longrightarrow 0$$
where each map is defined by the external product :
$$u\longmapsto u\wedge(f_1e_1+f_2e_2+...+f_re_r)$$
Hence, in degree $T$, the following is an exact sequence of vector spaces :
\begin{align*}
  &0\longrightarrow C_{T-t_1-t_2-...-t_r}\longrightarrow\bigoplus_{i=1}^rC_{T-t_1-...-\widehat{t_i}-...-t_r}\longrightarrow\cdots\\
  &\cdots\longrightarrow\bigoplus_{i<j}C_{T-t_i-t_j}\longrightarrow\bigoplus_{i=1}^rC_{T-t_i}\longrightarrow C_T\longrightarrow (C/(f_1,f_2,...,f_r))_T\longrightarrow 0
\end{align*}
From this one could calculate dimensions of the vector spaces involved:
$$\dim (C/(f_1,f_2,...,f_r))_T=\dim C_T-\sum_{i=1}^r\dim C_{T-t_i}+\sum_{i<j}\dim C_{T-t_i-t_j}-...$$
where we recognize the alternate sum already known to Bézout and Cayley.

\paragraph{Jouanolou} Hurwitz's project encompasses those of his predecessors Bézout, Hesse, Sylvester, Cayley, Mertens: one should describe every homogeneous component of the ideal of \textit{Trägheitsformen} of $(f_1,...,f_r)$, and, if possible, give a basis for each. These homogeneous components are $K[a,b,...]$-modules, and it is thus, essentially, linear algebra. Eventhough, Hurwitz concedes that ``the problem of determining all \textit{Trägheitsformen} of a module [=ideal] presents important difficulties''\footnote{\textit{Cf.} \cite{hurwitz1913}, p.~614.}. This project has been tackled again by Jean-Pierre Jouanolou in 1980.

Jouanolou studied the Koszul complex of the ideal $(f_1,...,f_r)$ in the language of Grothendieck's theory of schemes. Here, the problem of the \textit{Trägheitsformen} becomes a problem of ``local cohomology''. Now, we shall only explain how one of the propositions above by Hurwitz translates into this conceptual frame. Jouanolou demonstrates\footnote{\textit{Cf.} \cite{jouanolou1980} \S~2.1 to 2.8.} that some groups of cohomology with support in the ideal $\mathfrak{M}=(X_0,...,X_n)$ are null; in order to do so, he uses spectral sequences abutting to the hypercohomology of Koszul complex relative to the fonctor $\Gamma_{\mathfrak{M}}$ of sections with support in $\mathfrak{M}$. Thus, for example, if $n>r$, one has $H_{\mathfrak{M}}^0(C/(f_1,...,f_r))=0$; in other words, the ideal $(f_1,...,f_n)$ has no proper \textit{Trägheitsform}.

\section{Toric varieties}\label{Toric}

A theorem by D. N. Bernshtein \cite{bernshtein1975} published in 1975 
also gives the degree of the eliminand for systems of equations with support
in a convex set. Moreover, this theorem leads to the same kind of alternate sum as was found in Bézout's calculations. Bernshtein uses tools and concepts unknown to Bézout (infinite series in several variables, Minkowski's volume). The following year, an article \cite{kushnirenko1976} by Kushnirenko shows how to build a Koszul complex that leads directly to the afore-mentioned alternate sum. It is closely related to the theory of ``toric varieties'' developed in the 1970's. 

We shall now use toric varieties and a Koszul complex ; but rather than giving a full account of Kushnirenko's results about equations with support in a convex set, we shall merely give a complete proof of Bézout's theorem for $n$ incomplete equations of the second species (\textit{cf.} section \ref{Demonstration2} above), thereby also subsiding to the gap in Bézout's own demonstration.

Let there be $r$ incomplete equations of the second species, with $n$ unknowns :
$$\left\lbrace\begin{array}{l}
f^{(1)}=0\\
f^{(2)}=0\\
\vdots\\
f^{(r)}=0
\end{array}\right.$$
such that $\text{supp}(f^{(i)})=E_{t^{(i)},a^{(i)},b^{(i)}}$ for $1\leq i\leq n$, using the notations on p.~\pageref{secondspecies} above. The theory of toric varieties will provide us with an algebraic variety $X(\Delta)$ which is a compactification of the torus $(\mathbb{C}^\times)^n$ obtained by glueing affine varieties. This variety will allow a geometric interpretation of the vector spaces of polynomials studied by Bézout, as sets of global sections of some fiber bundles on $X(\Delta)$. Let us start with a few preliminaries.

\paragraph{Proposition 1} Let $P$ be the convex envelop in $\mathbb{R}^n$ of the support $E_{t,a,b}$ of an incomplete equation of the second species. Then $E_{t,a,b}=P\cap\mathbb{Z}^n$.

\emph{Demonstration.} $E_{t,a,b}$ is defined on p.~\pageref{secondspecies} by a system of inequations describing a convex polytope $P'$ in $\mathbb{R}^n$. If $P$ is the convex envelop in $\mathbb{R}^n$ of $E_{t,a,b}$, one must have $P\subset P'$. To prove equality, it is enough to check that every vertex of $P'$ belongs to $E_{t,a,b}$. By saturating some of the inequations defining $P'$, one can easily calculate its vertices. It has $(n^2+2n-3)$ vertices, each vertex belonging to $n$ faces of dimension 1. Some of those vertices may coincide in degenerate cases, \textit{i. e.} when $t$, $a$, $b$ verify certain relations. There are nine classes of vertices :
\begin{enumerate}[(i)]
\item $(0,0,...,0)$
\item $(a_1,b-a_1,0,0,...,0)$
\item $(b-a_2,a_2,0,0,...,0)$
\item for every $1\leq i\leq n$, one vertex $(0,...,0,a_i,0,...,0)$
\item for every $3\leq i\leq n$, one vertex $(a_1,b-a_1,0,...,0,t-b,0,...,0)$ where $t-b$ is the $i$-th coordinate
\item for every $3\leq i\leq n$, one vertex $(b-a_2,a_2,0,...,0,t-b,0,...,0)$ where $t-b$ is the $i$-th coordinate
\item for every $3\leq i\leq n$, one vertex $(a_1,0,...,0,t-a_1,0,...,0)$ where $t-a_1$ is the $i$-th coordinate
\item for every $3\leq i\leq n$, one vertex $(0,a_2,0,...,0,t-a_2,0,...,0)$ where $t-a_2$ is the $i$-th coordinate
\item for every $3\leq i\leq n$, and $1\leq j\leq n$ with $j\neq i$, one vertex $(0,...,0,a_i,0,...,0,t-a_i,0,...,0)$ where $a_i$ is the $i$-th coordinate and $(t-a_i)$ is the $j$-th coordinate
\end{enumerate}
Each of those vertices clearly belong to $\mathbb{Z}^n$ hence to $E_{t,a,b}$. \textit{q. e. d.}

\paragraph{Proposition 2} Let $P$ and $\Pi$ be the convex envelops in $\mathbb{R}^n$ of the supports $E_{t,a,b}$ and $E_{\theta,\alpha,\beta}$ of two incomplete equations of the second species. Then $E_{t+\theta,a+\alpha,b+\beta}$ is also the support of an incomplete equation of the second species, and the polytope $(P+\Pi)$ is its convex envelop in $\mathbb{R}^n$.

\emph{Demonstration.} By linearity, it is clear that $(t+\theta,a+\alpha,b+\beta)$ verify the ``restrictive conditions'' on p.~\pageref{secondspecies}. Let $P'$ be the convex envelop of $E_{t+\theta,a+\alpha,b+\beta}$ in $\mathbb{R}^n$. The polytope $(P+\Pi)$ is defined by
$$P+\Pi=\lbrace y\in\mathbb{R}^n\mid (\exists x\in P,\,\xi\in\Pi)\ y=x+\xi\rbrace$$
Now using the calculations in the demonstration of prop. 1, one sees\footnote{Beware that this crucial fact won't work for the third species of incomplete equations.} that all vertices of $P'$ belong to $P+\Pi$. Hence $P'\subset P+\Pi$.

On the other hand, as $P$, $\Pi$ and $P'$ are all defined by sets of \emph{linear} inequations similar to those on p.~\pageref{secondspecies}, it is also obvious that every $x+\xi\in P+\Pi$ belongs to $P'$. Hence $P'=P+\Pi$. \textit{q. e. d.}

\paragraph{Example} For $n=3$, such polytopes have 8 faces and 12 vertices, \textit{cf.} an example on figure 5.

\paragraph{Construction of $X(\Delta)$} We briefly recall the construction of an algebraic variety over $\mathbb{C}$ associated with a fan. The theory of toric varieties associates to every strongly rational convex cone $\sigma$ an affine algebraic variety $U_\sigma$. A strongly rational convex cone is a subset $\sigma$ of $\mathbb{R}^n$ generated over $\mathbb{R}_+$ by a finite family of vectors with rational coordinates, and such that $\sigma\cap(-\sigma)=\lbrace 0\rbrace$. A fan is a family $\Delta$ of such cones, such that each face of a cone in $\Delta$ is also a cone in $\Delta$, and that the intersection of two cones in $\Delta$ is a face of each. A fan $\Delta$ gives rise to an algebraic variety $X(\Delta)$ by glueing together the corresponding affine pieces. We are going to work with a fan $\Delta$ closely related to the polytopes described above.

\paragraph{Description of the fan $\Delta$} Let $\lbrace e_1,e_2,...,e_n\rbrace$ be the canonical basis in $\mathbb{R}^n$. The maximal cones of the fan $\Delta$ are simplicial cones generated over $\mathbb{R}_+$ by families of $n$ vectors. There are $(n^2+2n-3)$ such cones, corresponding to the following families of $n$ vectors :
\begin{itemize}
\item the cone generated by family $\lbrace e_1,e_2,...,e_n\rbrace$
\item the cone generated by family $\lbrace -e_1,-e_1-e_2,e_3,e_4,...,e_n\rbrace$
\item the cone generated by family $\lbrace -e_2,-e_1-e_2,e_3,e_4,...,e_n\rbrace$
\item the $n$ cones generated by families of the form
  $$\lbrace e_1,e_2,...,\widehat{e_i},...,e_n\rbrace\cup\lbrace -e_i\rbrace$$
\item the $2(n-2)$ cones generated by families of the form
  $$\lbrace e_3,e_4,...,\widehat{e_i},...,e_n\rbrace\cup\lbrace -e_1,-e_1-e_2,-e_1-e_2-...-e_n\rbrace$$
  or of the form
  $$\lbrace e_3,e_4,...,\widehat{e_i},...,e_n\rbrace\cup\lbrace -e_2,-e_1-e_2,-e_1-e_2-...-e_n\rbrace$$
\item the $2(n-2)$ cones generated by families of the form
  $$\lbrace e_3,e_4,...,\widehat{e_i},...,e_n\rbrace\cup\lbrace -e_1,e_2,-e_1-e_2-...-e_n\rbrace$$
  or of the form
  $$\lbrace e_3,e_4,...,\widehat{e_i},...,e_n\rbrace\cup\lbrace -e_2,e_1,-e_1-e_2-...-e_n\rbrace$$
\item the $(n-2)(n-1)$ cones generated by families of the form
  $$(\lbrace e_1,e_2,...,e_n\rbrace-\lbrace e_i,e_j\rbrace)\cup\lbrace -e_i,-e_1-e_2-...-e_n\rbrace$$
  where $i\geq 3$ and $j\neq i$.
\end{itemize}
The fan $\Delta$ is the set of all cones generated by sub-families of these families.

\paragraph{Remark} This fan could also be described as the fan of cones over the faces of a polytope dual to the polytope $P$ occuring in the previous propositions, and it has been calculated as such. See \cite{fulton1993} p.~26. As a matter of fact, the resulting fan does not depend upon the particular choice of $t,a,b$. There is a correspondance $\sigma\mapsto u(\sigma)$ between the maximal cones of $\Delta$ and the vertices of $P$ : \label{vertices}
\begin{itemize}
\item if $\sigma$ is generated by $\lbrace e_1,e_2,...,e_n\rbrace$, then $u(\sigma)=(0,0,...,0)$
\item if $\sigma$ is generated by $\lbrace -e_1,-e_1-e_2,e_3,e_4,...,e_n\rbrace$, then $u(\sigma)=(a_1,b-a_1,0,0,...,0)$
\item if $\sigma$ is generated by $\lbrace -e_2,-e_1-e_2,e_3,e_4,...,e_n\rbrace$, then $u(\sigma)=(b-a_2,a_2,0,0,...,0)$
\item if $\sigma$ is generated by $\lbrace e_1,e_2,...,\widehat{e_i},...,e_n\rbrace\cup\lbrace -e_i\rbrace$, then $u(\sigma)=(0,...,0,a_i,0,...,0)$
\item if $\sigma$ is generated by $\lbrace e_3,e_4,...,\widehat{e_i},...,e_n\rbrace\cup\lbrace -e_1,-e_1-e_2,-e_1-e_2-...-e_n\rbrace$, then $u(\sigma)=(a_1,b-a_1,0,...,0,t-b,0,...,0)$ ($t-b$ is the $i$-th coordinate)
\item if $\sigma$ is generated by $\lbrace e_3,e_4,...,\widehat{e_i},...,e_n\rbrace\cup\lbrace -e_2,-e_1-e_2,-e_1-e_2-...-e_n\rbrace$, then $u(\sigma)=(b-a_2,a_2,0,...,0,t-b,0,...,0)$
\item if $\sigma$ is generated by $\lbrace e_3,e_4,...,\widehat{e_i},...,e_n\rbrace\cup\lbrace -e_1,e_2,-e_1-e_2-...-e_n\rbrace$ then $u(\sigma)=(a_1,0,...,0,t-a_1,0,...,0)$
\item if $\sigma$ is generated by $\lbrace e_3,e_4,...,\widehat{e_i},...,e_n\rbrace\cup\lbrace -e_2,e_1,-e_1-e_2-...-e_n\rbrace$ then $u(\sigma)=(0,a_2,0,...,0,t-a_2,0,...,0)$
\item if $\sigma$ is generated by $(\lbrace e_1,e_2,...,e_n\rbrace-\lbrace e_i,e_j\rbrace)\cup\lbrace -e_i,-e_1-e_2-...-e_n\rbrace$ where $i\geq 3$ and $j\neq i$, then $u(\sigma)=(0,...,0,a_i,0,...,0,t-a_i,0,...,0)$ ($a_i$ is the $i$-th coordinate, and $t-a_i$ is the $j$-th coordinate)
\end{itemize}
When considering several polytopes $P^{(i)}$, we shall write $u^{(i)}(\sigma)$ the vertex of $P^{(i)}$ corresponding to the maximal cone $\sigma$. For a polytope $\Pi$, we shall use the greek letter $\upsilon(\sigma)$.

\paragraph{Example} For $n=3$, see a representation of the fan $\Delta$ on figure 6: each triangle represents one maximal cone.

\paragraph{Affine open sets $U_\sigma\subset X(\Delta)$} For each cone $\sigma$ of $\Delta$, one puts
$$U_\sigma=\text{Spec}(A_\sigma)$$
where $A_\sigma\subset\mathbb{C}\lbrack\chi_1,\chi_1^{-1},\chi_2,\chi_2^{-1},...,\chi_n,\chi_n^{-1}\rbrack$ is defined by :
$$A_\sigma=\bigoplus_{\begin{array}{c}\scriptstyle u\in\mathbb{Z}^n,\\[-.1cm]\scriptstyle(\forall v\in\sigma)\ \left<u,v\right>\geq 0\end{array}}\mathbb{C}\chi^u$$
Here $\chi_1,\chi_2,...,\chi_n$ are the indeterminates over the base field $\mathbb{C}$. If $\tau\in\Delta$ is a face of $\sigma\in\Delta$, there is a
natural mapping
$U_\tau\rightarrow U_\sigma$ embedding $U_\tau$ as a principal open subset of $U_\sigma$ (\textit{cf.} \cite{fulton1993}, p.~18), and one can thus build a variety $X(\Delta)$ by glueing all affine pieces $U_{\sigma_1}$ and $U_{\sigma_2}$ along $U_{\sigma_1}\cap U_{\sigma_2}=U_{\sigma_1\cap\sigma_2}$. In other words :
$$X(\Delta)=\lim_{\sigma\in\Delta}\text{ind.}(\text{Spec}(A_\sigma))$$

\paragraph{Construction of a vector bundle $O(D_P)$ on $X(\Delta)$} If $P$ is the convex envelop in $\mathbb{R}^n$ of the support of an incomplete equation of the second species, one can define a line bundle $O(D_P)$ on $X(\Delta)$ as follows. The bundle is trivial on each $U_\sigma$, \textit{ie.} $\simeq\mathbb{C}\times U_\sigma$. On the intersection of two maximal cones $\sigma_1$ and $\sigma_2$, the change of map is given by :
$$
\begin{diagram}
  \node{\mathbb{C}\times U_{\sigma_1}} \node{\mathbb{C}\times U_{\sigma_2}}\\[2]
  \node{\mathbb{C}\times U_{\sigma_1\cap\sigma_2}} \arrow[2]{n,J} \arrow{e} \node{\mathbb{C}\times U_{\sigma_1\cap\sigma_2}} \arrow[2]{n,J}\\
  \node{(t,\,x)} \arrow{e,T} \node{(\chi^{u(\sigma_1)-u(\sigma_2)}(x)t,\,x)}
\end{diagram}
$$
These are isomorphisms because $\chi^{u(\sigma_1)-u(\sigma_2)}$ is a unit of $A_{\sigma_1\cap\sigma_2}$. The compatibility of the changes of map on $U_{\sigma_1\cap\sigma_2\cap\sigma_3}$ comes from the fact that $\chi^{u(\sigma_1)-u(\sigma_3)}=\chi^{u(\sigma_2)-u(\sigma_3)}\chi^{u(\sigma_1)-u(\sigma_2)}$.

In other words, the sheaf of germs of sections of the vector bundle is isomorphic to the ideal sheaf generated over $A_\sigma$ by $\chi^{u(\sigma)}\in\mathbb{C}[\chi_1,\chi_1^{-1},...,\chi_n,\chi_n^{-1}]$ for every maximal cone $\sigma\in\Delta$.

\paragraph{Proposition 3} Under the same conditions as proposition 1, one has
$$\Gamma(X(\Delta),O(D_P))=\bigoplus_{u\in P\cap\mathbb{Z}^n}\mathbb{C}\chi^u$$

\emph{Demonstration.} On one hand, we must prove that $\chi^u$ is a regular section of $O(D_P)$ over $U_\sigma$, \textit{i. e.} $\chi^{u-u(\sigma)}\in A_\sigma$, for every maximal cone $\sigma$ and every $u\in P\cap\mathbb{Z}^n$. It is easily verified when $u$ is a vertex of $P$ using the description of the vertices on p.~\pageref{vertices}, and this is enough. On the other hand, if $u\in\mathbb{Z}^n$ verifies $(\forall\sigma)\ \chi^{u-u(\sigma)}\in A_\sigma$, one must prove that $u\in P$. Suppose it is not the case. Hahn-Banach theorem implies the existence of a hyperplane separating $u$ from the convex polytope $P$ : there exists $v\in\mathbb{R}^n$ and $a\in\mathbb{R}$ such that
$$(\forall u'\in P)\ \left<u',v\right> > a\quad\text{but}\quad \left<u,v\right> < a.$$
For the cone $\sigma$ containing $v$, this contradicts $(\forall\sigma)\ \chi^{u-u(\sigma)}\in A_\sigma$. Hence it is impossible that $u\not\in P$.

\paragraph{Proposition 4} Let $P^{(1)}$ and $P^{(2)}$ be the convex envelops of the supports of two incomplete equations of the second species, with $P^{(1)}\cap\mathbb{Z}^n=E_{t^{(1)},a^{(1)},b^{(1)}}$ and $P^{(2)}\cap\mathbb{Z}^n=E_{t^{(2)},a^{(2)},b^{(2)}}$. Let $f^{(2)}\in\Gamma(X(\Delta),O(D_{P^{(2)}}))$. Multiplication by $f^{(2)}$ induces a natural map
$$\begin{diagram}\node{O(D_{P^{(1)}})}\arrow{e,t}{\times f^{(2)}}\node{O(D_{P^{(1)}+P^{(2)}})}\end{diagram}$$

\emph{Demonstration.} This map is defined locally on each affine subset $U_\sigma$ by :
$$
\begin{diagram}
  \node{\phi\in\Gamma(U_\sigma,O(D_{P^{(1)}}))} \arrow{e,T} \arrow{s,T} \node{\phi\chi^{-u^{(1)}(\sigma)}\in A_\sigma} \arrow{s,T}\\
    \node{\phi f^{(2)}\in\Gamma(U_\sigma,O(D_{P^{(1)}+P^{(2)}}))} \arrow{e,T} \node{\phi\chi^{-u^{(1)}(\sigma)}f^{(2)}\chi^{-u^{(2)}(\sigma)}\in A_\sigma}
\end{diagram}
$$
In the following, we shall write $\widetilde{f}^{(2)}=f\chi^{-u^{(2)}(\sigma)}$.

\paragraph{Koszul complex} One could thus define, locally on every open affine subset $U_\sigma$, a whole complex isomorphic to a Koszul complex. Let $\Pi$ be the convex envelop of the support of any incomplete equation of the second species\footnote{Avoid degenerate cases where some of the vertices $\upsilon(\sigma)$ coincide.}, and $f^{(1)}$, $f^{(2)}$,..., $f^{(r)}$ as above. Let $L_\sigma$ be the $A_\sigma$-module defined by
$$L_\sigma=\bigoplus_{i=1}^rA_\sigma$$
The Koszul complex gives a sequence of maps of $A_\sigma$-modules :
$$
\begin{diagram}
  \node{0} \arrow{s} \node{0} \arrow{s}\\
  \node{\Gamma(U_\sigma,O(D_\Pi))} \arrow{s} \arrow{e,t}{\simeq} \node{\Lambda^0L_\sigma=A_\sigma} \arrow{s,r}{\wedge\widetilde{f}}\\
  \node{\Gamma(U_\sigma,\bigoplus_{i=1}^rO(D_{\Pi+P^{(i)}}))} \arrow{s} \arrow{e,t}{\simeq} \node{\Lambda^1L_\sigma=L_\sigma} \arrow{s,r}{\wedge\widetilde{f}}\\
  \node{\Gamma(U_\sigma,\bigoplus_{i<j}O(D_{\Pi+P^{(i)}+P^{(j)}}))} \arrow{s} \arrow{e,t}{\simeq} \node{\Lambda^2L_\sigma} \arrow{s,r}{\wedge\widetilde{f}}\\
  \node{\vdots} \arrow{s} \node{\vdots} \arrow{s,r}{\wedge\widetilde{f}}\\
  \node{\Gamma(U_\sigma,\bigoplus_{i=1}^rO(D_{\Pi+P^{(1)}+...+\widehat{P^{(i)}}+...+P^{(r)}}))} \arrow{s} \arrow{e,t}{\simeq} \node{\Lambda^{r-1}L_\sigma} \arrow{s,r}{\wedge\widetilde{f}}\\
  \node{\Gamma(U_\sigma,O(D_{\Pi+P^{(1)}+...+P^{(r)}}))} \arrow{e,t}{\simeq} \node{\Lambda^rL_\sigma\simeq A_\sigma}
\end{diagram}
$$
where $\widetilde{f}=\begin{pmatrix} f^{(1)}\chi^{-u^{(1)}(\sigma)}\\
f^{(2)}\chi^{-u^{(2)}(\sigma)}\\
\vdots\\
f^{(n)}\chi^{-u^{(n)}(\sigma)}
\end{pmatrix}$. Those maps glue with each other on the intersections $U_{\sigma_1}\cap U_{\sigma_2}$ into a sequence of maps of sheaves over $X(\Delta)$:
$$0 \arrow{e} O(D_\Pi) \arrow{e} \bigoplus_{i=1}^rO(D_{\Pi+P^{(i)}})) \arrow{e} \bigoplus_{i<j}O(D_{\Pi+P^{(i)}+P^{(j)}})) \arrow{e} \hdots \arrow{e} O(D_{\Pi+P^{(1)}+...+P^{(r)}})$$

\paragraph{Theorem 1} This sequence of maps of sheaves over $X(\Delta)$ is an exact sequence.

\emph{Demonstration}. We prove it locally on every open affine subset $U_\sigma$. In fact, each $\widetilde{f}^{(i)}$ has a non-zero constant term because $\chi^{u^{(i)}}$ belongs to the support of $f^{(i)}$. One can thus use the following trick by Mertens, in order to prove that $\widetilde{f}^{(1)},\widetilde{f}^{(2)},...,\widetilde{f}^{(r)}$ is a regular sequence. One must prove that, for all $s\leq r$, $\widetilde{f}^{(s)}$ does not divide zero modulo $(\widetilde{f}^{(1)},\widetilde{f}^{(2)},...,\widetilde{f}^{(s-1)})$. Suppose
$$\sum_{i=1}^s\phi^{(i)}\widetilde{f}^{(i)}=0$$
Let us recall that our base field is an extension of $\mathbb{Q}$, and that the coefficients of our polynomials $f^{(i)}$ are indeterminates over $\mathbb{Q}$ (\textit{cf.} p.~\pageref{secondspecies}). The constant term in $\widetilde{f}^{(i)}$ is such an indeterminate, call it $c^{(i)}$. The base field $K$ could thus be written as $k(c^{(1)},...,c^{(r)})$ where $k$ is an extension of $\mathbb{Q}$. There is an isomorphism\footnote{This is inspired by Mertens \cite{mertens1886} p.~528-529.}:
$$\begin{diagram}
  \node{k\left[c^{(1)},...,c^{(r)}\right]\left[\chi_1,\chi_1^{-1},...,\chi_n,\chi_n^{-1}\right]/\left(\widetilde{f}^{(1)},...,\widetilde{f}^{(s-1)}\right)} \arrow{s,l}{\simeq} \arrow{e,!} \node{c^{(i)}} \arrow{s,T}\\
      \node{k\left[c^{(s)},...,c^{(r)}\right]\left[\chi_1,\chi_1^{-1},...,\chi_n,\chi_n^{-1}\right]} \arrow{e,!} \node{\left\lbrace\begin{array}{l}c^{(i)}-\widetilde{f}^{(i)}\text{ if }i<s\\
        c^{(i)}\text{ if }i\geq s\end{array}\right.}
\end{diagram}$$
The bottom ring is a polynomial ring, it is an integral domain, and thus
$$\phi^{(s)}\widetilde{f}^{(s)}\in(\widetilde{f}^{(1)},...,\widetilde{f}^{(s-1)})\quad\Rightarrow\quad\phi^{(s)}\in(\widetilde{f}^{(1)},...,\widetilde{f}^{(s-1)})$$
\textit{q. e. d.}

\paragraph{Theorem 2} For large enough $k$, the fibre bundle $O(D_{k\Pi})$ is very ample. As a consequence, $X(\Delta)$ is a projective variety embedded in  $\mathbb{P}^{\vert k\Pi\cap\mathbb{Z}^n\vert-1}(\mathbb{C})$. For such $k$, there exists $N$ such that the following sequence is exact :
$$0 \arrow{e} \Gamma(X,O(D_{(Nk+1)\Pi})) \arrow{e} \bigoplus_{i=1}^r\Gamma(X,O(D_{(Nk+1)\Pi+P^{(i)}})) \arrow{e} \bigoplus_{i<j}\Gamma(X,O(D_{(Nk+1)\Pi+P^{(i)}+P^{(j)}}))...$$

\emph{Demonstration}. The existence of $k$ and of a very ample $O(D_{k\Pi})$ is a consequence of the non-degeneracy of $\Pi$. See \cite{fulton1993} p.~69-70 for a proof. For such $k$, write $O(1)=O(D_{k\Pi})$. A theorem of Serre states that, if $\mathcal{F}$ is a coherent algebraic sheaf on a projective variety, then for $N$ large enough, $\mathcal{F}\otimes O(N)$ has trivial cohomology. This implies that, for $N$ large enough, the following sequence is exact :
$$0 \arrow{e} \Gamma(X,O(D_{\Pi})\otimes O(N)) \arrow{e} \bigoplus_i\Gamma(X,O(D_{\Pi+P^{(i)}})\otimes O(N))) \arrow{e}...$$
(\textit{cf.} Serre's theorem in \cite{grothendieck1963} 2.2.1, and its corollary 2.2.3).

\paragraph{Corollary (Bézout's theorem for the second species)} For $k$ and $N$ as in theorem 2, when $n=r$, the dimension of the cokernel of the last map\footnote{\textit{i. e.} the map defined by $\Lambda^{n-1}L_\sigma\rightarrow\Lambda^nL_\sigma$ over every $U_\sigma$.} is
$$\prod_{i=1}^nt^{(i)}-\sum_{j=1}^n\prod_{i=1}^n(t^{(i)}-a_j^{(i)})+\prod_{i=1}^n(t^{(i)}-b^{(i)})-\sum_{i=1}^n\left\lbrack(a_1^{(i)}+a_2^{(i)}-b^{(i)})\prod_{j\neq i}(t^{(j)}-b^{(j)})\right\rbrack$$
It is an upper bound on the degree of the eliminand of $f^{(1)},...,f^{(n)}$.

\emph{Demonstration}. The alternate sum of dimensions of the vector spaces
appearing in the exact sequence above can be expressed as the following finite difference of order $n$:
$$\Delta_{t^{(n)},a^{(n)},b^{(n)}}...\Delta_{t^{(2)},a^{(2)},b^{(2)}}\Delta_{t^{(1)},a^{(1)},b^{(1)}}\vert E_{T,A,B}\vert$$
where $T,A,B$ are the degrees occuring in $\Lambda^nM$, \textit{i. e.}:
\begin{align*}
  T&=(Nk+1)\theta+t^{(1)}+t^{(2)}+...+t^{(n)}\\
  (\forall i)\quad A_i&=(Nk+1)\alpha_i+a_i^{(1)}+a_i^{(2)}+...+a_i^{(n)}\\
  B&=(Nk+1)\beta+b^{(1)}+b^{(2)}+...+b^{(n)}
\end{align*}
Calculating $\vert E_{T,A,B}\vert$ is an easy combinatorial problem. One has:
\begin{align*}
  \vert E_{T,A,B}\vert=&{T+n \choose n}-\sum_{i=1}^n{T-A_i+n-1 \choose n}\\
  &+{T-B+n-2 \choose n}-(A_1+A_2-B){T-B+n-2 \choose n-1}
\end{align*}
Let us recall the well-known formula
$${m\choose p}-{m-n\choose p}=\sum_{k=1}^n{m-k\choose p-1}$$
as well as the finite difference of a product:
$$\Delta_t(P(T)Q(T))=Q(T)\Delta_tP(T)+P(T-t)\Delta_tQ(T)$$
Using these formulae enables us to calculate the quantity above and prove the result stated.

\paragraph{Proof of the statement on p.~\pageref{statement}} For $r$ polynomials with $n$ indeterminates, the map $(f^{(1)},f^{(2)},...f^{(r)})_{\leq T,A,B}$ of section \ref{statement} is none other than the last map\footnote{\textit{i. e.} the map defined by $\Lambda^{r-1}L_\sigma\rightarrow\Lambda^rL_\sigma$ over every $U_\sigma$.} of the sequence in theorem~2. For values of $T$, $A$, $B$ ensuring the non-degeneracy of $\Pi$, and for $k$ and $N$ as in theorem~2, the sequence is exact, so the kernel of this map must be the image of the next map : if $(\phi^{(1)},...,\phi^{(r)})$ is an element of the kernel, then there exists a family of polynomials
$$(\psi^{(ij)})_{i,j}\in\bigoplus_{i<j}C_{\leq T-t^{(i)}-t^{(j)},A-a^{(i)}-a^{(j)},B-b^{(i)}-b^{(j)}}$$
such that, in particular,
$$\phi^{(1)}=\sum_{j=2}^r(-1)^j\psi^{(1j)}f^{(j)}$$
\textit{q. e. d.}

\paragraph{Third species of incomplete equations, $n=r=3$} For polynomials $f^{(1)},f^{(2)},f^{(3)}$ of the third species in three indeterminates, most of the arguments above are still valid, although there is a major problem with proposition 2. In the demonstration of proposition 2, the coordinates of the vertices of $P'$ were linear in $t+\theta,a+\alpha,b+\beta$ and could be written as the sums of the coordinates of the corresponding vertices of $P$ and $\Pi$, the three polytopes having the same form; but, as we said in section \ref{Demonstration3} p.~\pageref{Demonstration3}, there are eight different forms of polynomials of the third species. The convex envelops of their supports are polytopes of differents forms (\textit{cf.} figure 4). In order to fix the demonstration of proposition 2, one is going to study a larger class of polytopes, from which the eight forms of polytopes of the third species are only degenerate forms. These polytopes are represented on figure 7. If $(t,a,b)$ belongs to the third species, such a polytope is the convex envelop of a set $E_{t,a,b,s}\subset\mathbb{Z}^3$ defined by :
$$\left\lbrace\begin{aligned}
&0\leq k_1\leq a_1,\quad 0\leq k_2\leq a_2,\quad 0\leq k_3\leq a_3,\\
&k_1+k_2\leq b_3,\quad k_1+k_3\leq b_2,\quad k_2+k_3\leq b_1,\\
&k_1+k_2+k_3\leq t,\\
&2k_1+k_2+k_3\leq s_1,\quad k_1+2k_2+k_3\leq s_2,\quad k_1+k_2+2k_3\leq s_3
\end{aligned}\right.$$

\paragraph{Proposition 5} If $(t,a,b)$ belongs to the third species, put
for $1\leq i\leq 3$:
$$s_i=\min(t+a_i,\ b_{i+1}+b_{i+2})$$
where the indices are modulo 3 (for example $b_4=b_1$). Then $E_{t,a,b,s}=E_{t,a,b}$.

\emph{Demonstration.} It is clear that $E_{t,a,b,s}\subset E_{t,a,b}$. Moreover, if $(k_1,k_2,k_3)\in E_{t,a,b}$, one has:
\begin{align*}
  2k_i+k_{i+1}+k_{i+2}=(k_1+k_2+k_3)+k_i\leq t+a_i\\
  2k_i+k_{i+1}+k_{i+2}=(k_i+k_{i+1})+(k_i+k_{i+2})\leq b_{i+2}+b_{i+1}
\end{align*}
Hence $2k_i+k_{i+1}+k_{i+2}\leq\min(t+a_i,b_{i+1}+b_{i+2})=s_i$, which concludes the proof.

\paragraph{A new fan} One subdivides the fan $\Delta$ on figure 6, using new rays through the following vectors:
$$-e_1-e_3,\quad -e_2-e_3,\quad -2e_1-e_2-e_3,\quad -e_1-2e_2-e_3,\quad -e_1-e_2-2e_3$$
The new fan is represented on figure 8; it is compatible with the new class of polytopes. Propositions 1 to 4 and theorem 1 and 2 are valid for this fan and these polytopes.

\paragraph{Proposition 6} If $P$ is the convex envelop of $E_{T,A,B,S}$, then
\begin{align*}\vert P\cap\mathbb{Z}^3\vert=&{T+3 \choose 3}+\sum_{i=1}^3\left\lbrack{T-B_i+1 \choose 3}-{T-A_i+2 \choose 3}\right\rbrack\\
  &-\sum_{i=1}^3\left\lbrack(A_{i+1}+A_{i+2}-B_i){T-B_i+1 \choose 2}\right\rbrack\\
  &+\sum_{i=1}^3\left\lbrack(T+A_i-B_{i+1}-B_{i+2}+1){T+A_i-S_i+1 \choose 2}-2{T+A_i-S_i+2 \choose 3}\right\rbrack
\end{align*}

\emph{Demonstration.} This combinatorial formula derives by truncation from any of the eight formulae given in section \ref{Demonstration3} for polytopes of the third species.

\paragraph{Remark} On the other way around, by specializing $S_i=\min(T+A_i,B_{i+1}+B_{i+2})$ in the formula above, one could also derive the eight formulae given in section \ref{Demonstration3}. For example, if $T+A_3>B_1+B_2$, the corresponding term in the expression above is :
\begin{align*}
  &(T+A_3-B_1-B_2+1){T+A_3-B_1-B_2+1 \choose 2}-2{T+A_3-B_1-B_2+2 \choose 3}\\
  &={T+A_3-B_1-B_2+1 \choose 3}
\end{align*}
as an easy computation would reveal. This term appears, as it should, in the formulae for the 2d, the 4th, the 6th and the 8th forms.

\paragraph{Corollary} Analog to theorem 2 and its corollary, when $k$ and $N$ are large enough, the dimension of the cokernel of the last map in the Koszul complex gives the following upper bound on the degree of the eliminand:
\begin{align*}&\prod_{i=1}^3t^{(i)}+\sum_{i=1}^3\left\lbrack\prod_{j=1}^3(t^{(j)}-b_i^{(j)})-\prod_{j=1}^3(t^{(j)}-a_i^{(j)})\right\rbrack\\
  &-\sum_{i=1}^3\sum_{j=1}^3(a_{i+1}^{(j)}+a_{i+2}^{(j)}-b_i^{(j)})\prod_{k\neq j}(t^{(k)}-b_i^{(k)})\\
  &+\sum_{i=1}^3\left\lbrack\sum_{j=1}^3(t^{(j)}+a_i^{(j)}-b_{i+1}^{(j)}-b_{i+2}^{(j)})\prod_{k\neq j}(t^{(k)}+a_i^{(j)}-s_i^{(j)})-2\prod_{j=1}^3(t^{(j)}+a_i^{(j)}-s_i^{(j)})\right\rbrack
\end{align*}
where $s_i^{(j)}=\min(t^{(j)}+a_i^{(j)},b_{i+1}^{(j)}+b_{i+2}^{(j)}).$

\emph{Demonstration.} Use prop. 6 and calculate:
$$\Delta_{t^{(3)},a^{(3)},b^{(3)},s^{(3)}}\Delta_{t^{(2)},a^{(2)},b^{(2)},s^{(2)}}\Delta_{t^{(1)},a^{(1)},b^{(1)},s^{(1)}}\vert E_{T,A,B,S}\vert.$$

\paragraph{Bézout's own formula for the third species} In his treatise, Bézout does not use the truncated polytopes that we have described above. He calculates everything under the hypothesis that all polytopes appearing in the exact sequence belong to one and the same form, among the eight forms pertaining to the third species of incomplete equations. He thus finds eight different formulae, that could be derived from the formula of previous corollary by specializing the $s_i^{(j)}$ to their corresponding values.

One might doubt that any of those eight formulae could apply to the cases where $f^{(1)}$, $f^{(2)}$ and $f^{(3)}$ belong to distinct forms, because the 9 parameters $s_i^{(j)}$ could each be specialized in two distinct ways (either $t^{(j)}+a_i^{(j)}$ or $b_{i+1}^{(j)}+b_{i+2}^{(j)}$), which makes 256 possible outcomes. Nevertheless, calculation reveals that, for every $i$, the last term in square brackets in the sum above, \textit{i. e.}
$$\sum_{j=1}^3(t^{(j)}+a_i^{(j)}-b_{i+1}^{(j)}-b_{i+2}^{(j)})\prod_{k\neq j}(t^{(k)}+a_i^{(j)}-s_i^{(j)})-2\prod_{j=1}^3(t^{(j)}+a_i^{(j)}-s_i^{(j)}),$$
only takes two possible values after such specialization. Indeed, write
$h_i^{(j)}=t^{(j)}+a_i^{(j)}-b_{i+1}^{(j)}-b_{i+2}^{(j)}$. For a given $i$, if $h_i^{(j)}\leq 0$ for two or three among the three possible values of the index $j$, the term in square brackets vanishes identically ; but if $h_i^{(j)}\leq 0$ for at most one value of $j$, then the term in square brackets is equal to $h_i^{(1)}h_i^{(2)}h_i^{(3)}$. Finally, the upper bound on the degree of the eliminand is:
\begin{align*}&\prod_{i=1}^3t^{(i)}+\sum_{i=1}^3\left\lbrack\prod_{j=1}^3(t^{(j)}-b_i^{(j)})-\prod_{j=1}^3(t^{(j)}-a_i^{(j)})\right\rbrack\\
  &-\sum_{i=1}^3\sum_{j=1}^3(a_{i+1}^{(j)}+a_{i+2}^{(j)}-b_i^{(j)})\prod_{k\neq j}(t^{(k)}-b_i^{(k)})+\sum_{i=1}^3\left\lbrack\varepsilon_i\prod_{j=1}^3h_i^{(j)}\right\rbrack
\end{align*}
where $\varepsilon_i=0$ or $1$. This coincides with the eight formulae given by Bézout in his treatise \cite{bezout1779} \S~119-127.

\section{Comparing methods}
The geometrical origin of Cayley's researches, unlike Bézout's, might obscure the identity of methods. Both scholars met with the same mechanisms of linear algebra, having to do with the same unsolved problem : the exactness of a sequence of linear maps. Could Cayley, in some way or another, have known of Bézout's treatise ? He does not mention it ; but he knew of Waring's \textit{Meditationes algebraicae}, the second edition of which contains a very brief summary of Bézout's ideas. Cayley's method is exactly the same as Bézout's, informed by the theory of determinants, Sylvester's dialytic method, and the new-born matrix symbolism. The alternate sum of dimensions is present in Bézout's, in Waring's, and in Cayley's works.

Despite of these similarities, with Sylvester, Hesse and Cayley, elimination theory is on a new track characterized by:
\begin{itemize}
\item The important role of projective and algebraic geometry in focusing on homogeneous polynomials.
\item The annexion of elimination theory within the growing \emph{theory of invariants} and the systematic search for invariants.
\item The calculatory trend aiming at explicit formulas, mainly depending on determinants and using matrix algebra.
\end{itemize}

It would be misleading to see the concept of ideal where it is not. Ideals appeared in algebra at the cross-road with number theory in a research trail starting with the \textit{Disquisitiones Arithmeticae} of Gauss, up to Kummer, Dedekind, Weber and Kronecker at the end of XIXth century. Yet, there is novelty in Bézout's treatment of ``sum-equations'' in 1779. This novelty and the lack of rigour in Bézout's treatise, as well as its refusal of geometry\footnote{In this respect, Euler had clearly seen the relation between elimination and projection on the axis of a cartesian coordinate system. The problem of particular cases due to points of intersection at infinity could not be solved before the introduction of projective methods in algebraic geometry.}, had endangered its reception. These obstacles partially overcome, Hesse's and Cayley's articles definitely give a posterity to Bézout's treatise.

The peculiar dialectic between general statements and the many generic cases was also present in Bézout's treatise, and it is both a weakness and a strength. For exemple, the fact that the degree of the eliminand is always less than or equal to the product of the degrees of the given equations, is a general statement. A perfectly grounded universal proof of this statement had to await until the end of XIXth century (Serret, Schmidt, Hurwitz); but Bézout had already understood that this upper bound is the exact degree of the eliminand in the generic case of $n$ ``complete equations''. He had also known of other generic cases where the exact degree of the eliminand is less than the product of the degrees and could be precisely acertained. As proven above, the formulae found by Bézout in many cases are confirmed by the theory of toric varieties and the method of Bernshtein and Kushnirenko.


\section{Appendix : an elementary proof for the first species of incomplete equations}
Let us consider a system of three incomplete equations of first species :
$$\left\lbrace\begin{array}{l}
f^{(1)}=0\\
f^{(2)}=0\\
f^{(3)}=0
\end{array}\right.$$
We use the same notations as above.

For any set of integers $T$, $A_1$, $A_2$, $A_3$ verifying the conditions pertaining to the first species of incomplete equations, let us call $(f)_{\leq T,A}$ the linear map defined by
$$\begin{array}{lrcl}
(f)_{\leq T,A}:&\displaystyle\bigoplus_{i=1}^3C_{\leq T-t^{(i)},A-a^{(i)}}&\longrightarrow &C_{\leq T,A}\\
&(\phi^{(1)},\phi^{(2)},\phi^{(3)})&\longmapsto &\sum_{i=1}^3\phi^{(i)}f^{(i)}
\end{array}$$

We are going to build a resolution of $\text{coker}(f)_{\leq T,A}$, \textit{i. e.} an exact sequence of linear maps :
$$\begin{array}{rl}
&0\\
&\downarrow\\
&C_{\leq T-t^{(1)}-t^{(2)}-t^{(3)},\ A-a^{(1)}-a^{(2)}-a^{(3)}}\\
(h)_{\leq T,A}&\downarrow\\
&\bigoplus\limits_{i=1}^3C_{\leq T-t^{(1)}-t^{(2)}-t^{(3)}+t^{(i)},\ A-a^{(1)}-a^{(2)}-a^{(3)}+a^{(i)}}\\
(g)_{\leq T,A}&\downarrow\\
&\bigoplus\limits_{i=1}^3C_{\leq T-t^{(i)}, A-a^{(i)}}\\
(f)_{\leq T,A}&\downarrow\\
&C_{\leq T,A}
\end{array}$$

\paragraph{The kernel of $(f)_{\leq T,A}$}
Suppose
$$\sum_{i=1}^3\phi^{(i)}f^{(i)}=0$$
There is an isomorphism of $\mathbb{Q}[x_1,x_2,x_3]$-algebras%
\footnote{This is inspired by Mertens \cite{mertens1886} p.~528-529.}
$$\begin{array}{rcl}
\mathbb{Q}\left[(u_{i,k})_{2\leq i\leq 3, k\in\text{supp}(f^{(i)})}\right]\left[x_1,x_2,x_3\right]/\left(f^{(2)},f^{(3)}\right)&\simeq &\mathbb{Q}\left[(u_{i,k})_{k\neq(0,0,0)}\right]\left[x_1,x_2,x_3\right]\\
u_{i,k}&\mapsto &\left\lbrace\begin{array}{l}u_{i,(0,0,0)}-f^{(i)}\text{ if }k=(0,0,0)\\
u_{i,k}\text{ if }k\neq(0,0,0)\end{array}\right.
\end{array}$$
The ring on the right-hand side is a polynomial ring, it is an integral domain, and
$$\phi^{(1)}f^{(1)}\in(f^{(2)},f^{(3)})\quad\Rightarrow\quad\phi^{(1)}\in(f^{(2)},f^{(3)})$$
Hence there exists $\psi^{(2)}$, $\psi^{(3)}$ such that
$$\phi^{(1)}=\psi^{(3)}f^{(2)}-\psi^{(2)}f^{(3)}$$
First of all, we are going to prove that it is possible to choose such $\psi^{(2)}$, $\psi^{(3)}$ with
$$\psi^{(2)}\in C_{\leq T-t^{(1)}-t^{(3)},A-a^{(1)}-a^{(3)}},\quad
\psi^{(3)}\in C_{\leq T-t^{(1)}-t^{(2)},A-a^{(1)}-a^{(2)}}$$
Indeed, if $\deg\psi^{(3)}>T-t^{(1)}-t^{(2)}$, then we should have
$$[\psi^{(3)}][f^{(2)}]-[\psi^{(2)}][f^{(3)}]=0$$
where the brackets designate the terms of highest total degree of a polynomial. Now, $f^{(2)}$ and $f^{(3)}$ being generic, the greatest common divisor of $[f^{(2)}]$ and $[f^{(3)}]$ over $K=\mathbb{Q}((u_{i,k})_{i,k})$ must be a monomial. But no monomial divides any of those two polynomials. So:
$$(\exists\lambda\in K)\quad[\psi^{(3)}]=\lambda[f^{(3)}],\quad[\psi^{(2)}]=\lambda[f^{(2)}]$$
We could then write
$$\phi^{(1)}=(\psi^{(3)}-\lambda f^{(3)})f^{(2)}-(\psi^{(2)}-\lambda f^{(2)})f^{(3)}$$
and the two new multiplier-polynomials are of total degree less than the total degree of the old ones. We thus prove by induction that there exists two multiplier-polynomials $\psi^{(2)}$, $\psi^{(3)}$ with
$$\deg\psi^{(2)}\leq T-t^{(1)}-t^{(3)},\quad\deg\psi^{(3)}\leq T-t^{(1)}-t^{(2)},$$
$$\phi^{(1)}=\psi^{(3)}f^{(2)}-\psi^{(2)}f^{(3)}$$

Let us now suppose that $\deg_1\psi^{(3)}>A_1-a_1^{(1)}-a_1^{(2)}$, then we should have
$$[\psi^{(3)}]_1[f^{(2)}]_1-[\psi^{(2)}]_1[f^{(3)}]_1=0$$
where the brackets $[\cdot]_1$ designate the terms of highest degree in $x_1$. Now, $f^{(2)}$ and $f^{(3)}$ being generic, one has:
$$[f^{(2)}]_1=x_1^{a^{(2)}_1}F^{(2)}$$
$$[f^{(3)}]_1=x_1^{a^{(3)}_1}F^{(3)}$$
where $F^{(2)}$ and $F^{(3)}$ are irreducible polynomials over $K$. Thus there exists a polynomial $\Lambda$ in $x_2$, $x_3$ such that:
$$[\psi^{(3)}]_1=\Lambda F^{(3)}x_1^{a_1^{(3)}-\min\left(a_1^{(2)},a_1^{(3)}\right)}$$
$$[\psi^{(2)}]_1=\Lambda F^{(2)}x_1^{a_1^{(2)}-\min\left(a_1^{(2)},a_1^{(3)}\right)}$$ 
Because of the hypothesis on $\deg_1\psi^{(3)}$, as soon as $A_1$ will be large enough, $\Lambda$ will be divisible by $x_1^{\min\left(a_1^{(2)},a_1^{(3)}\right)}$, and thus one will be able to write
$$\phi^{(1)}=\left(\psi^{(3)}-\frac{\Lambda}{x_1^{\min\left(a_1^{(2)},a_1^{(3)}\right)}}f^{(3)}\right)f^{(2)}-\left(\psi^{(2)}-\frac{\Lambda}{x_1^{\min\left(a_1^{(2)},a_1^{(3)}\right)}}f^{(2)}\right)f^{(3)}$$
where the two multiplier-polynomials are of degree in $x_1$ less than the old ones. We thus recursively prove that there exists two multiplier-polynomials $\psi^{(2)}$, $\psi^{(3)}$ with
$$\deg_1\psi^{(2)}\leq A_1-a_1^{(1)}-a_1^{(3)},\quad\deg_1\psi^{(3)}\leq A_1-a_1^{(1)}-a^{(2)}$$
We could have worked in the same fashion with $\deg_2$ and with $\deg_3$. Their should still be one care: let us make sure that the transformation of the multiplier-polynomials described above in order to decrease their degrees with respect to a single unknown, let us say $x_1$, does not increase their degrees with respect with to $x_2$ or $x_3$. Concerning $x_2$ (the same reasoning holds for $x_3$), one has:
$$\begin{array}{rcl}
\deg_2(\Lambda f^{(3)})&=&\deg_2[\psi^{(3)}]_1-\deg_2F^{(3)}+a_2^{(3)}\\
&\leq& (\deg\psi^{(3)}-\deg_1\psi^{(3)})-(t^{(3)}-a_1^{(3)})+a_2^{(3)}\\
&<& (T-t^{(1)}-t^{(2)})-(A_1-a_1^{(1)}-a_1^{(2)})-(t^{(3)}-a_1^{(3)})+a_2^{(3)}
\end{array}$$
If $T-t^{(1)}-t^{(2)}-t^{(3)}\leq (A_1-a_1^{(1)}-a_1^{(2)}-a_1^{(3)})+(A_2-a_2^{(1)}-a_2^{(2)}-a_2^{(3)})$, then one has, as we wished:
$$\deg_2(\Lambda f^{(3)})< A_2-a_2^{(1)}-a_2^{(2)}$$
We have thus achieved the proof that it is possible to choose $\psi^{(2)}$, $\psi^{(3)}$ with
$$\psi^{(2)}\in C_{\leq T-t^{(1)}-t^{(3)},A-a^{(1)}-a^{(3)}},\quad
\psi^{(3)}\in C_{\leq T-t^{(1)}-t^{(2)},A-a^{(1)}-a^{(2)}}$$
$$\phi^{(1)}=\psi^{(3)}f^{(2)}-\psi^{(2)}f^{(3)}$$
Let there be such $\phi^{(1)}=\psi^{(3)}f^{(2)}-\psi^{(2)}f^{(3)}$. Obviously
$$(\phi^{(1)},-\psi^{(3)}f^{(1)},\psi^{(2)}f^{(1)})\in\ker(f)_{\leq T,A}$$
The other elements of $\ker(f)_{\leq T,A}$ with first coordinate equal to $\phi^{(1)}$ can all be written
$$(\phi^{(1)},\phi^{(2)}-\psi^{(3)}f^{(1)},\phi^{(3)}+\psi^{(2)}f^{(1)})$$
where $(0,\phi^{(2)},\phi^{(3)})\in\ker(f)_{\leq T,A}$, that is to say
$$\phi^{(2)}f^{(2)}+\phi^{(3)}f^{(3)}=0$$
As $f^{(2)}$ and $f^{(3)}$ are irreducible polynomials over $K$, then there exists $\psi^{(1)}$ such that
$$\phi^{(2)}=\psi^{(1)}f^{(3)},\quad\phi^{(3)}=-\psi^{(1)}f^{(2)}$$
Hence the kernel of $(f)_{\leq T,A}$ is the image of the following linear map $(g)_{\leq T,A}$:
$$\begin{array}{rrcl}
(g)_{\leq T,A}:&\bigoplus\limits_{i=1}^3C_{\leq T-t^{(1)}-t^{(2)}-t^{(3)}+t^{(i)},\ A-a^{(1)}-a^{(2)}-a^{(3)}+a^{(i)}}&\longrightarrow&\bigoplus\limits_{i=1}^3C_{\leq T-t^{(i)},A-a^{(i)}}\\
&(\psi^{(1)},\psi^{(2)},\psi^{(3)})&\longmapsto&(\psi^{(3)}f^{(2)}-\psi^{(2)}f^{(3)},\\
&&&\psi^{(1)}f^{(3)}-\psi^{(3)}f^{(1)},\\
&&&\psi^{(2)}f^{(1)}-\psi^{(1)}f^{(2)})
\end{array}$$

\paragraph{The kernel of $(g)_{\leq T,A}$}
Let $(\psi^{(1)},\psi^{(2)},\psi^{(3)})\in\ker(g)_{T,A}$. The polynomials $f^{(i)}$ are irreducible, so that, again:
$$(\exists\Lambda\in C_{\leq T-t^{(1)}-t^{(2)}-t^{(3)},\ A-a^{(1)}-a^{(2)}-a^{(3)}})\quad\psi^{(1)}=\Lambda f^{(1)},\ \psi^{(2)}=\Lambda f^{(2)},\ \psi^{(3)}=\Lambda f^{(3)}$$
Hence the kernel of $(g)_{\leq T,A}$ is the image of the following linear map $(h)_{\leq T,A}$:
$$\begin{array}{rrcl}
(h)_{\leq T,A}:&C_{\leq T-t^{(1)}-t^{(2)}-t^{(3)},\ A-a^{(1)}-a^{(2)}-a^{(3)}}&\longrightarrow&\bigoplus\limits_{i=1}^3C_{\leq T-t^{(1)}-t^{(2)}-t^{(3)}+t^{(i)},\ A-a^{(1)}-a^{(2)}-a^{(3)}+a^{(i)}}\\
&\Lambda&\longmapsto&(\Lambda\psi^{(1)},\Lambda\psi^{(2)},\Lambda\psi^{(3)})
\end{array}$$
Moreover, this map is injective. In other words, we now have an exact sequence of linear maps building a resolution of $\text{coker}(f)_{\leq T,A}$:
$$\begin{array}{rl}
&0\\
&\downarrow\\
&C_{\leq T-t^{(1)}-t^{(2)}-t^{(3)},\ A-a^{(1)}-a^{(2)}-a^{(3)}}\\
(h)_{\leq T,A}&\downarrow\\
&\bigoplus\limits_{i=1}^3C_{\leq T-t^{(1)}-t^{(2)}-t^{(3)}+t^{(i)},\ A-a^{(1)}-a^{(2)}-a^{(3)}+a^{(i)}}\\
(g)_{\leq T,A}&\downarrow\\
&\bigoplus\limits_{i=1}^3C_{\leq T-t^{(i)}, A-a^{(i)}}\\
(f)_{\leq T,A}&\downarrow\\
&C_{\leq T,A}
\end{array}$$

\paragraph{Finite differences and alternate sum}\label{difference}
Using several times the rank theorem, one thus has:
$$\begin{array}{rcl}\dim\text{coker}(f)_{\leq T,A}&=&\dim C_{\leq T,A}\\
&&-\dim\bigoplus\limits_{i=1}^3C_{\leq T-t^{(i)},A-a^{(i)}}\\
&&+\dim\bigoplus\limits_{i=1}^3C_{\leq T-t^{(1)}-t^{(2)}-t^{(3)}+t^{(i)},\ A-a^{(1)}-a^{(2)}-a^{(3)}+a^{(i)}}\\
&&-\dim C_{\leq T-t^{(1)}-t^{(2)}-t^{(3)},\ A-a^{(1)}-a^{(2)}-a^{(3)}}
\end{array}$$
This expression could be rewritten into a single finite difference of order 3:
$$\begin{array}{rcl}
\dim\text{coker}(f)_{\leq T,A}&=&((\dim C_{\leq T,A}-\dim C_{\leq T-t^{(1)},A-a^{(1)}})\\
&&-(\dim C_{\leq T-t^{(2)},A-a^{(2)}}-\dim C_{\leq T-t^{(1)}-t^{(2)},A-a^{(1)}-a^{(2)}}))\\
&&-((\dim C_{\leq T-t^{(3)},A-a^{(3)}}-\dim C_{\leq T-t^{(1)}-t^{(3)},A-a^{(1)}-a^{(3)}})\\
&&-(\dim C_{\leq T-t^{(2)}-t^{(3)},A-t^{(2)}-t^{(3)}}-\dim C_{\leq T-t^{(1)}-t^{(2)}-t^{(3)},A-a^{(1)}-a^{(2)}-a^{(3)}}))\\
&=&\Delta_{t^{(3)},a^{(3)}}\Delta_{t^{(2)},a^{(2)}}\Delta_{t^{(1)},a^{(1)}}\dim C_{\leq T,A}
\end{array}$$
Calculating $\dim C_{\leq T,A}$ is an easy combinatorial problem:
$$\dim C_{\leq T,A}={T+3 \choose 3}+{T-A_1+2 \choose 3}+{T-A_2+2 \choose 3}+{T-A_3+2 \choose 3}$$
Let us recall a well-known formula:
$${m\choose p}-{m-n\choose p}=\sum_{k=1}^n{m-k\choose p-1}$$
Using this formula several times enables one to calculate the finite difference above. It is a constant, independant of $T$ and $A$. Eventually,
$$\dim\text{coker}(f)_{T,A}=t^{(1)}t^{(2)}t^{(3)}-\sum_{i=1}^3(t^{(1)}-a^{(1)}_i)(t^{(2)}-a^{(2)}_i)(t^{(3)}-a^{(3)}_i)$$
If there exists an eliminand in $x_1$ of lowest degree, that number is an upper bound on its degree.

\includepdf[pages=-]{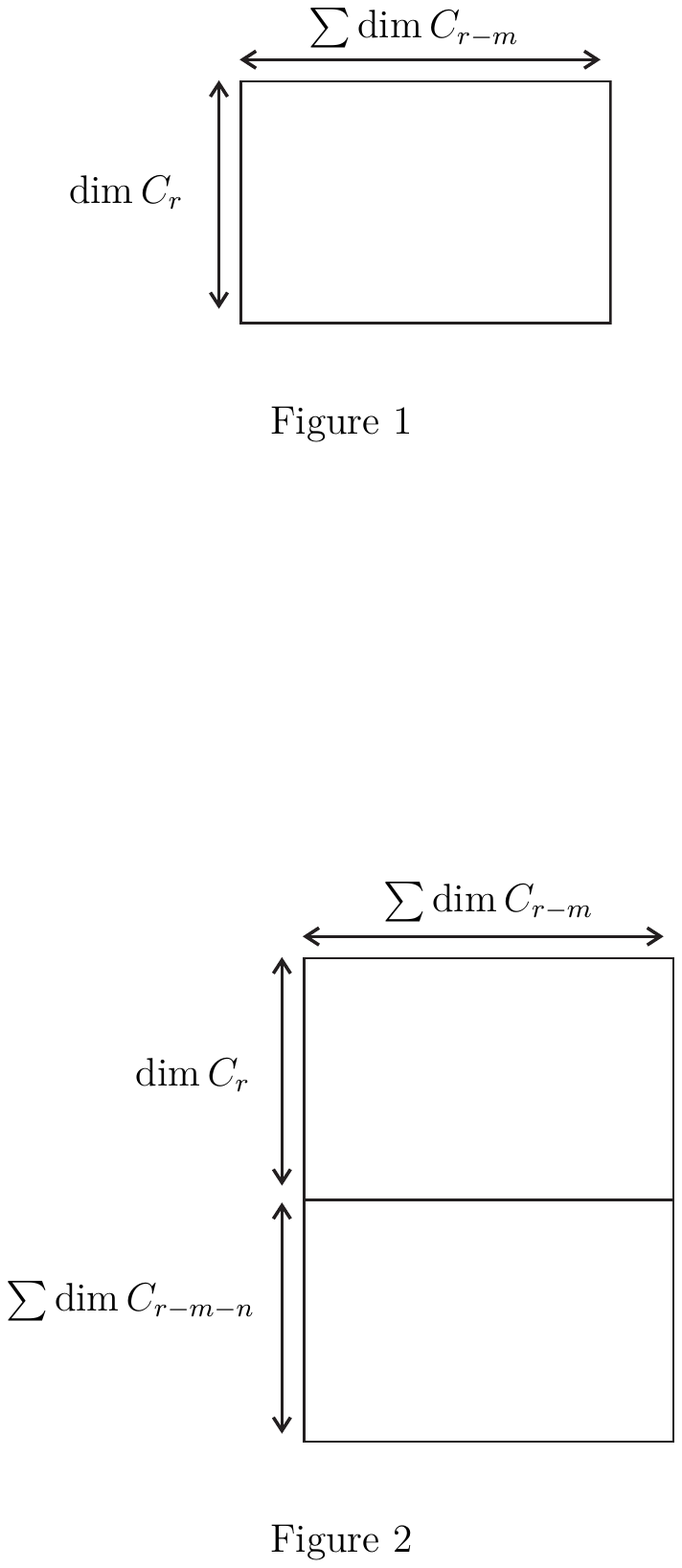}

\end{document}